\documentclass{article}
\usepackage{amssymb}
\usepackage{amsmath}
\usepackage{setspace}
\newtheorem{theorem}{Theorem}

\newtheorem{lemma}[theorem]{Lemma}

\newtheorem{remark}[theorem]{Remark}

\begin{document}

\begin{center}
{\LARGE At the origins and in the vanguard of peri-dynamics, non-local }

{\LARGE and higher gradient continuum mechanics. }

{\LARGE An underestimated and still topical contribution of Gabrio Piola}

by Francesco dell'Isola, Ugo Andreaus and Luca Placidi

\bigskip
\end{center}

Gabrio Piola's scientific papers have been underestimated in the
mathematical-physics literature. Indeed a careful reading of them proves
that they are original, deep and far reaching. Actually -even if his
contribution to mechanical sciences is not completely ignored- one can
undoubtedly say that the greatest part of his novel contributions to
mechanics, although having provided a great impetus and substantial
influence on the work of many preminent mechanicians, is in fact generally
ignored. It has to be remarked that authors \cite{capecchiruta2007}
dedicated many efforts to the aim of unveiling the true value of Gabrio
Piola as a scientist; however, some deep parts of his scientific results
remain not yet sufficiently illustrated. Our aim is to prove that non-local
and higher gradient continuum mechanics was conceived already in Piola's
works and to try to explain the reasons of the unfortunate circumstance
which caused the erasure of the memory of this aspect of Piola's
contribution. Some relevant differential relationships obtained in Piola 
\cite{piola1845-6} are carefully discussed, as they are still nowadays too
often ignored in the continuum mechanics literature and can be considered
the starting point of Levi-Civita's theory of Connection for Riemannian
manifolds.

\section{Introduction}

Piola's contribution to the mechanical sciences is not completely ignored:
indeed his contribution to the formulation of balance equations of force in
Lagrangian description is universally recognized (when and how this first
re-discovery of Piola occurred will be the object of another investigation).
In this context the spirit of Piola's works can be recognized in many modern
contributions (see e.g. \cite{quiligottietal2003}). One can undoubtedly say
that the greatest part of his novel contributions to Mechanics, although
having imparted a great momentum to and substantial influence on the work of
many prominent mechanicians, is in fact generally ignored.

Although the last statement may seem at first sight exaggerated, the aim of
the present paper is exactly to prove it while presenting the evidence of a
circumstance which may seem surprising: some parts the works of Gabrio Piola
represent a topical contribution as late as the year 2013.

Those who have appreciated the works of Russo \cite{russo2013}, \cite%
{russo2003} will not be at all shocked by such a statement, as there is
evidence that many scientific contributions remained unsurpassed for
centuries, if not millennia. Therefore one thesis that we want to put
forward in this paper is that the contribution of Gabrio Piola should not be
studied with the attitude of the historian of science but rather with the
mathematical rigor needed to understand a contemporary textbook or a
research paper.

On the other hand, the authors question the concept of "historical method"
especially when applied to history of science and history of mathematics. We
claim that there is not any peculiar "historical method" to be distinguished
from the generic "scientific method" which has to be applied to describe any
other kind of phenomena, although the subject of the investigation is as
complex as those involved in the transmission, storage and advancement of
scientific knowledge. A fortiori, however, imagine that one could determine
precisely in what constitutes such a "historical method": then it MUST
include the capability of the historician to understand, master and
reconstruct rigorously the mathematical theories which he has decided to
study from the historical point of view. In other words: a historician of a
particular branch of mathematics has to master completely the theory whose
historical development he wants to describe. It is rather impossible, for
instance, that somebody who does not know the theory of integration could
recognize that (see \cite{russo2003}) Archimedes actually used rigorous
arguments leading to the proof of the existence of the integral of a
quadratic function. Moreover, together with the linguistic barriers (one has
to know doric Greek to understand Archimedes and XIXth century Italian to
understand Gabrio Piola), there are also notational difficulties: one should
not naively believe that HIS\ OWN notations are advanced and modern while
the notations found in the sources are clumsy and primitive. Actually
notations are a matter of "arbitrary choice" and from this point of view -
remember that mathematics is based on axiomatic definition of abstract
concepts to which the mathematician assigns a meaning by means of axioms and
definitions - all notations are equally acceptable. Very often historians of
mathematics \footnote{%
See for instance Russo (2003) page 53 and ff. for what concerns the
difficulties found by historicians who did not know trigonometry in
recognizing that Hellenist science had formulated it but with different
fundamental variables and notations} decide that a theory is much more
modern than it actually is, simply because they do not find "the modern
symbols" or the "modern nomenclature" in old textbooks. For instance, if one
does not find in a textbook the symbol $\int $ , this does not mean that the
integral was not known to the author of that textbook. It could simply mean
that the technology of printing at the age of that textbook required the use
of another symbol or of another symbolic method. Indeed some formulas by
Lagrange or Piola seem at first sight to the authors of the present paper to
resemble, for their complicated length, lines of commands for LaTeX.
Actually the historian has to READ carefully the textbooks which he wants to
assess and interpret: when these books are books whose content is a
mathematical theory, reading them implies reading all the fundamental
definitions, lemmas and properties, which are needed to follow its logical
development.

In the authors' opinion, in Truesdell and Toupin \cite{truesdelltoupin1960}
the contribution of Piola to mechanical sciences is accounted for only
partially while in Truesdell \cite{truesdell1968} it is simply overlooked.
It has to be remarked that authors \cite{capecchiruta2007} dedicated many
efforts to the aim of unveiling the true value of Gabrio Piola as a
scientist; however, some deep parts of his scientific results remain not yet
sufficiently illustrated.

Our aim is

\begin{itemize}
\item to prove that non-local and higher gradients continuum mechanics was
conceived already in Piola's works starting from a clever use of the
Principle of Virtual Work

\item to explain the unfortunate circumstances which caused the erasure from
memory of this aspect of Piola's contribution, although his pupils respected
so greatly his scientific standing that they managed to dedicate an
important square to his name in Milan (close to the Politecnico) to be named
after him, while a statue celebrating him was erected in the Brera Palace,
also in Milan.
\end{itemize}

Finally some differential relationships obtained in \cite{piola1845-6} are
carefully discussed, as they are still nowadays too often used without
proper attribution in the continuum mechanics literature and can be
considered the starting point of the Levi-Civita theory of Connection in
Riemannian manifolds.

The main source for the present paper is the work\medskip

\bfseries

Piola, G., Memoria intorno alle equazioni fondamentali del movimento di
corpi qualsivogliono considerati secondo la naturale loro forma e
costituzione, Modena, Tipi del R.D. Camera, (1845-1846),

\mdseries\medskip

\noindent but the authors have also consulted other works by Piola \cite%
{piola1825}, \cite{piola1833}, \cite{piola1835}, \cite{piola1856}.

In all the above-cited papers by Piola the kinematical descriptor used is
simply the placement field defined on the reference configuration: in these
works there is no trace of more generalized models of the type introduced by
the Cosserats \cite{cosseratcosserat1909}. However the spirit of Piola's
variational formulation (see, e.g., \cite{bedford1985}, \cite%
{berdichevsky2009}, \cite{dahermaugin1986}, \cite{germain1973a}, \cite%
{germain1973b}, \cite{kroner1968}, \cite{maugin2011}, \cite{toupin1962}, 
\cite{toupin1964}) and his methods for introducing generalized stress
tensors can be found in the papers by Green and Rivlin (\cite{greenrivlin1964a}, \cite{greenrivlin1965} \cite{greenrivlin1964c} and \cite{greenrivlin1964b}) and
also many modern works authored for instance by Neff and his co-workers, 
\cite{neff2006a}, \cite{neff2006b}, \cite{neffjeong2009} and of by Forest
and his co-workers \cite{forest2009}, \cite{forestetal2011}.

\section{Linguistic, ideological and cultural barriers impeding the
transmission of knowledge}

It is evident to many authors and it is very often recognized in the
scientific literature that linguistic barriers may play a negative role in
the transmission and advancement of science. We recall here, for instance,
that Peano \cite{peano1903} in 1903, being aware of the serious consequences
which a Babel effect can have on the effective collaboration among
scientistis, tried to push the scientific community towards the use of Latin
or of an especially constructed "lingua franca" in scientific literature,
i.e. the so-called latino sine flexione.

Actually in Russo \cite{russo2013}, \cite{russo2003} the author clearly
analyses the consequences of the existence of those linguistic ideological
and cultural barriers which did not permit the Latin speaking scholars to
understand the depth of Hellenistic science: the beginning of the economical
and social processes leading to the Middle Ages.

\subsection{Gabrio Piola as a protagonist of Italian Risorgimento
(Resurgence)}

It is surprising that some important contributions to mechanics of a
well-known scientist remained unnoticed and have been neglected for a so
long a time. Actually, after a careful observation of distinct traces and by
gathering hints and evidences, one can propose a well-founded conjecture:
Gabrio Piola has been a leading cultural and scientific protagonist of
Italian Risorgimento (Resurgence).

The main evidence of this statement has been found, e.g., in his eulogy in
memoriam of his "Maestro" Vincenzo Brunacci. This eulogy was written in
1818 (three years after the famous Rimini Proclamation by Giocchino Murat
that, to give an idea of its content, started with "Italians! The hour has
come to engage in your highest destiny" and which is generally considered as
the beginning of the Italian Resurgence). In this eulogy (completely
translated in the Appendix B) there is a continuous reference to the Italian
Nation which, in that time, could pursue some serious legal difficulties for
the author of such a eulogy, leading eventually to the loss of his
personal freedom.

The eulogy starts with the words

`` It is extremely painful for us to announce in
this document the death of a truly great man, who, as during his life, was a
glory for Italy, ''

and ends with the words

``May these last achievements of such an inventive
and ingenious Geometer be delivered up to a capable and educated scholar,
who could enlighten them as they deserve, for the advancement of SCIENCES,
for the glory of the AUTHOR\ and for the prestige of ITALY''.

In the body of the eulogy one can find the following statements:

\begin{itemize}
\item `` It seemed as if the Spirit of Italy who
was in great sufferance because in that time the most brilliant star of all
mathematical sciences, the illustrious Lagrangia, had left the Nation, that
Spirit wanted to have the rise of another star, which being born on the
banks of the river Arno, was bound to become the successor of the first one.
This consideration is presenting itself even more spontaneous by when we
will remark that Brunacci was the first admirer in Italy of the luminous
Lagrangian doctrines, the scientist who diffused and supported them, the
scientist who in his studies was always a very creative innovator in their
applications. His first Maestri were two famous Italians, Father Canovai and
the great geometer Pietro Paoli ''
\end{itemize}

We remark that here Piola refers to Lagrange by his true and original
Italian name, Lagrangia, that he refers to Italy as a unique cultural
entity, that he deemed to exist the ``Spirit of
Italy'', that he refers to as Italians two
scientists who were Professors in Pisa (outside the Kingdom of
Lombardy--Venetia, where Piola lived and worked) .

\begin{itemize}
\item ``It is not licit for me neglecting to
indicate another subject in which -with honored efforts- our professor
distinguished himself. The Journal of Physical Chemistry of Pavia was
illustrated in many of his pages by his erudite pen; I will content myself
to indicate here three Memoirs where he examines the doctrine of capillary
attraction of Monsieur Laplace, comparing it with that of Pessuti and where,
with his usual frankness which is originated by his being persuaded of how
well-founded was his case, he proves with his firm reasonings, \textbf{%
whatever it is said\ by the French geometers}, some propositions which are
of great praise for the mentioned Italian geometer.''
\end{itemize}

For the purposes of this paper, we note that Piola remarked that Brunacci
gave in these memoirs of the Journal of Physical Chemistry, the role of the
champion of Italian science to the Italian Pesutti as counterposed to the
french geometer M. Laplace.

Brunacci greatly influenced Piola's scientific formation and rigorously
cultivated his ingenious spirit, as Piola himself recognized in many places
of his works. Piola was initiated by Brunacci to Mathematical Analysis but
was immediately attracted -since his first original creations- to
Mathematical Physics, which he based on the Principle of Virtual Velocities
(as Lagrange called what has been later called the Principle of Virtual Work)

Actually the aim of the whole scientific activity of Gabrio Piola has been
to demonstrate that such a Principle can be considered the basis of the
Postulation of every Mechanical Theory, see e.g. the papers \cite%
{carcaterra2005}, \cite{carcaterraakay2011}, \cite{carcaterraakay2007}, \cite%
{carcaterraetal2006}, \cite{carcaterraetal2000}, \cite%
{carcaterrasestieri1995}, \cite{dellisolavidoli1998a}, \cite%
{dellisolavidoli1998}, \cite{maurinietal2004a}, \cite{maurinietal2004b} for
the inclusion of the dissipative effects. Indeed he developed -by using the
Lagrangian Postulation- modern continuum mechanics, being -to our knowledge-
the first author who introduced the dual in power of the gradient of
velocity in the referential description of a continuous body. The
coefficients of what will be recognized to be a distribution in the modern
sense (as defined by Schwartz) were to be identified later, after the
revolutionary theories introduced by Ricci and Levi-Civita, as a double
tensor, the \textbf{Piola stress tensor}.

Some of the results presented in Piola's works (e.g. those concerning
continua the strain energy of which depends on higher gradients of the
strain measures) can be regarded even nowadays as among the most advanced
available in the literature.

\subsection{\textbf{Piola's works did not receive their due attention
because they are written in Italian.}}

It is clear that the strongest limiting factor to the full recognition of
Piola's contribution to Continuum Mechanics must be found in his
"ideological" choice: the use of the Italian language. Moreover he used a
very elegant and erudite style which can be understood and appreciated only
by a few specialists and he did not care if his works would be translated
into other languages (as later was decided by Levi-Civita who -instead-
cared to have some of his works translated into English and who wrote
directly some others in French (see the works by Ricci-Curbastro and
Levi-Civita \cite{levicivita1929}, \cite%
{levicivita1931}, \cite{levicivita1925}, \cite{riccilevicivita1900}).

It is clear that, in a historical period when all scientists of a given
Nation were using their own language in higher studies, when in every
University the official spoken language was the National one and where all
textbooks, essays and scientific Memoirs were written in the mother language
of the authors, Piola could not accept to admit the inferiority of his own
mother language and decided to use it for publishing his works.

A well-founded conjecture about this linguistic choice can be advanced:
although Piola was surely fluent in French (he edited in Italian some works
by Cauchy and cites long French excerpts by Poisson) he decided ("per la
gloria dell'Italia") for the glory of Italy to use his mother language, in
an historical climate in which the Italian Nation was not yet the united and
independent and therefore was not able to self-determine its destiny. This
was a patriotic choice which was repaid by a nearly complete neglect of his
contribution to mechanical science, exacerbated by the fact that Italian
authors seem to have underestimated his contributions (for a detailed
analysis of this point see \cite{capecchiruta2007}).

From a general point of view, the linguistic barriers often play a very
puzzling role in the diffusion of ideas and theories. As discussed in \cite%
{russo2013}, \cite{russo2003} the diffusion of Hellenistic science actually
was slowed by the great barrier represented by the ignorance of the language
used, but not stopped. The information slowly flowed from East to West, and
althought it needed some centuries, and in the end, maybe translated into a
Latin difficult to understand, still keeping Greek nomenclature and
terminology, this science managed to pollinate the Italian and European
Renaissance; however, the linguistic transfer corresponded to a nearly
complete loss of the knowledge about the identity of the scientists who had
first formulated the ideas at the basis of the scientific revolution. Even
the true period of the appearance of the scientific method was finally
postponed for more than a millennium.

It has not to be considered astonishing, then, that the contribution of
Piola still is permeating the modern Continuum Mechanic literature, but is
generally misunderstood, also by those who know better his contributions.

Indeed linguistic barriers are very often insurmountable.

\subsection{The mathematics used by Piola in his mechanics treatises}

The mathematics used by Piola is in every aspect modern, except in a very
important point. Indeed as Levi-Civita's absolute calculus was invented many
years later, Piola's presentation proceeds firmly and rigorously but
encumbered by a very heavy component-wise notation, which in the eyes of a
modern mechanician conveys an undeserved appearance of primitiveness. The
reader should not believe that Piola would refuse (as some mechanicians
still do!) to use the powerful tools given to us by Levi-Civita. Indeed
-again as proven in the eulogy he wrote for honoring his "Maestro"
Vincenzo Brunacci- Gabrio Piola knew how important are the choice of the
right notation and the conceptual tools for the advancement of science.
Furthermore, he calls "obscurantists" those who refused the nominalistic and
conceptual improvements introduced by Lagrange in Mathematical Analysis.

Unfortunately Piola did not have available to him the tool he needed to
progress more quickly in his research. It is astonishing to discover how
many results he managed to obtain notwithstanding this limitation.

\section{Non-Local Continuum Theories in Piola's works}

In the work by Piola \cite{piola1845-6} the homogenized theory which is
deduced by means of the identification of powers in the discrete micro-model
and in the continuous macro-model can be called (in the language used by
Eringen \cite{eringen1999}, \cite{eringenedelen1972}) a non-local theory.
Also some Italian authors (see e.g. \cite{polizzotto2001}), who contributed
to the field with important papers, seem not to give explicit recognition
that they were reformulating (and extending) the results already found by
Piola.

In Appendix A we translate those parts of Piola's work which are most
relevant in the present context and in this section we translate into modern
symbols the formulas which the reader may find in such an appendix in their
original form. Moreover, we will recall in a less suggestive, but more
direct and modern, language the statements made by Piola.

It is our opinion that some of Piola's arguments can compete in depth and
generality, even nowadays, with those which can be found in some of the most
advanced modern presentations. Postponing the analysis of Piola's
homogeneization process to a subsequent investigation, we limit ourselves
here to describe the continuum model which he deduces from the Principle of
Virtual Velocities for a discrete mechanical system constituted by a finite
set of molecules, which he considers to be (or, because of his controversy
with Poisson, he must accept as) the most fundamental Principle in his
Postulation process.

In Piola \cite{piola1845-6} (Capo I, pag. 8) the reference configuration of
the considered deformable body is introduced by labelling each material
particle with the three Cartesian coordinates $\left( a,b,c\right) .$ It is
suggestive to remark that the same notation is used in Hellinger \cite%
{hellinger1913}, see e.g., pag.605. We will denote by the symbol $X$ the
position occupied by each of the considered material particles in the
reference configuration. The placement of the body is then described by the
set of three scalar functions (Capo I, pag.8 and then pages 11-14) 
\begin{equation*}
x\left( a,b,c\right) ,y\left( a,b,c\right) ,z\left( a,b,c\right)
\end{equation*}%
which, by using a compact notation, we will denote with the symbol $\chi $
mapping any point in the reference configuration into its position in the
actual one.

\subsection{Piola's non-local internal interactions}

In Capo VI, on page 149 Piola introduces:

``the quantity $\rho $ (equations (3),(5), (6))
has the value given by the equation

\begin{align}
\rho^{2} & =\left[ x\left( a+f,b+g,c+k\right) -x\left( a,b,c\right) \right]
^{2}  \notag \\
& +\left[ y\left( a+f,b+g,c+k\right) -y\left( a,b,c\right) \right] ^{2} 
\tag{8}  \label{8} \\
& +\left[ z\left( a+f,b+g,c+k\right) -z\left( a,b,c\right) \right] ^{2}. 
\notag
\end{align}

''

So by denoting with the symbol $\bar{X}$ the particle labelled by Piola with
the coordinates $(a+f$, $b+g$, $c+k)$ we have, in modern notation, that 
\begin{equation}
\rho ^{2}(X,\bar{X})=\left\Vert \chi (\bar{X})-\chi (X)\right\Vert ^{2}. 
\tag{8bis}  \label{DELTARHO2}
\end{equation}%
In Capo VI on page 150 we read the following expression for the internal
work, relative to a virtual displacement $\delta \chi ,$ followed by a very
clear remark:

``

\begin{equation}
\int da\int db\int dc\int df\int dg\int dk\cdot\frac{1}{2}K\delta\rho 
\tag{10}
\end{equation}

[...] In it the integration limits for the variables $f,g,k$ will depend on
the surfaces which bound the body in the antecedent configuration, and also
on the position of the molecule $m$, which is kept constant, that is by the
variables $a,b,c$ which after the first three will also vary in the same
domain.''

Here the scalar quantity $K$ is introduced as the \textit{intensity} of the
force (see the page 147 of the translated part of \cite{piola1845-6} in
Appendix A) exerted by the particle $\bar{X}$ on the particle $X$ and the 
$\frac12$
is present as the action reaction principle holds. The quantity $K$ is
assumed to depend on $\bar{X},X$ and $\rho $ and manifestly it is measured
in $\left[ N\left( m\right) ^{-6}\right] $ (SI Units). In the number 72
starting on page 150 (translated completely in the appendix A) Piola
discusses the physical meaning of this scalar quantity and consequently
establishes some restrictions on the constitutive equations which have to be
assigned to it. Indeed he refrains from any effort to obtain for it an
expression in terms of microscopic quantities and limits himself to require
its objectivity by assuming its dependence on $\rho$, an assumption which
will have in the sequel some important consequences. Moreover he argues that
if one wants to deal with continua more general than fluids (for a
discussion of this point one can have a look on the recent paper \cite%
{auffreyetal2014}) then it may depend (in a symmetric way) also on the
Lagrangian coordinates of both $\bar{X}$ and $X:$ therefore%
\begin{equation*}
K(\bar{X},X,\rho )=K(X,\bar{X},\rho ).
\end{equation*}

On Page 151,152 we then read some statements which cannot be rendered
clearer:

``As a consequence of what was were said up to now
we can, by adding up the two integrals (1), (10), and by replacing the
obtained sum in the first two parts of the general equation (1) num$%
{{}^\circ}%
$.16., formulate the equation which includes the whole molecular mechanics.
Before doing so we will remark that it is convenient to introduce the
following definition

\begin{equation}
\Lambda =\frac{1}{4}\frac{K}{\rho }  \tag{11}  \label{CAPOVI-11ENG2}
\end{equation}%
by means of which it will be possible to introduce the quantity $\Lambda
\delta \rho ^{2}$ instead of the quantity $\frac{1}{2}K\delta \rho $ in the
sextuple integral (10); and that inside this sextuple integral it is
suitable to isolate the part relative to the triple integral relative to the
variables $f,g,k,$ placing it under the same sign of triple integral with
respect to the variables $a,b,c$ \ which includes the first part of the
equation: which is manifestly allowed. In this way the aforementioned
general equation becomes

\begin{equation*}
\int da\int db\int dc\cdot\left\{ \left( X-\frac{d^{2}x}{dt^{2}}\right)
\delta x+\left( Y-\frac{d^{2}y}{dt^{2}}\right) \delta y+\left( Z-\frac{d^{2}z%
}{dt^{2}}\right) \delta z\right.
\end{equation*}

\begin{equation}
\left. +\int df\int dg\int dk\cdot \Lambda \delta \rho ^{2}\right\} +\Omega
=0  \tag{12}  \label{CAPOVI-12ENG2}
\end{equation}%
where it is intended that (as mentioned at the beginning of the num$%
{{}^\circ}%
$.71.) it is included in the $\Omega $ the whole part which may be
introduced because of the forces applied to surfaces, lines or
well-determined points and also because of particular conditions which may
oblige some points to belong to some given curve or surface. ''

Piola is aware of the technical difficulty which he could be obliged to
confront in order to calculate the first variation of a square root: as he
knows that these difficulties have no physical counterparts, instead of $K $
he introduces another constitutive quantity $\Lambda$ which is the dual in
work of the variation $\delta\rho^{2}.$

\begin{remark}
\label{BoundednessandAttenuation}\textbf{Boundedness and attenuation
assumptions on }$K$\textbf{\ and }$\Lambda.$ Note that Piola explicitly
assumes the summability of the function $\Lambda\delta\rho^{2}=\frac{1}{4}%
\frac {K}{\rho}\delta\rho^{2}=\frac{1}{2}K\delta\rho$ and the boundedness of
the function $K.$ As a consequence when $\rho$ is increasing then $\Lambda$
decreases.
\end{remark}

\begin{remark}
\textbf{Objectivity of Virtual Work.} Note that $\delta \rho ^{2}$ and $%
\Lambda (X,\bar{X},\rho )$ are invariant (see \cite{steigmann2003}) under
any change of observer and as Piola had repeadedly remarked, see e.g. Capo
IV, num.48, page 86-87, the expression for virtual work has to verify this
condition. Remark also that, as the work is a scalar, in this point Piola's
reasoning is rendered difficult by his ignorance of Levi-Civita's tensor
calculus. In another formalism the previous formula can be written as
follows 
\begin{equation}
\int_{\mathcal{B}}[\left( b_{m}(X)-a(X)\right) \delta \chi (X)+\left( \int_{%
\mathcal{B}}\Lambda (X,\bar{X},\rho )\delta \rho ^{2}\mu (\bar{X})d\bar{X}%
\right) ]\mu (X)dX+\delta W(\partial \mathcal{B})=0  \tag{12bis}
\label{WORKBALANCE}
\end{equation}%
where $\mathcal{B}$ is the considered body, $\partial \mathcal{B}$ its
boundary, $\mu $ is the volume mass density, $b_{m}(X)$ is the (volumic)
mass specific externally applied density of force, $a(X)$ the acceleration
of material point $X,$ and $\delta W(\partial \mathcal{B})$ the work
expended on the virtual displacement by actions on the boundary $\partial 
\mathcal{B}$ and eventually the first variations of the equations expressing
the applied constraints on that boundary times the corresponding Lagrange
multipliers.
\end{remark}

In Eringen \cite{eringenedelen1972}, \cite{eringen1999}, \cite{eringen2002},
the non-local continuum mechanics is founded on a Postulation based on
Principles of balance of mass, linear and angular momenta, energy and
entropy. However in \cite{eringen2002} a chapter on variational principles
is presented.

One can easily recognize by comparing, for example, the presentation in \cite%
{eringen2002} with (\ref{WORKBALANCE}) that in the works by Piola the
functional

\begin{equation}
\left( \int_{\mathcal{B}}\Lambda (X,\bar{X},\rho )\delta \rho ^{2}\mu (\bar{X%
})d\bar{X}\right)  \tag{N1}  \label{INTERNAL WORKATX}
\end{equation}%
is assumed to satisfy a slightly generalized version of what in \cite%
{eringen2002} pag. 34 is called the

\begin{center}
\textit{Smooth Neighborhood Hypothesis }
\end{center}

which reads (in Eringen's work the symbol $V$ is used with the same meaning
as our symbol $\mathcal{B}$, $X^{\prime}\mathcal{\ }$instead of $\bar{X},$ $%
x $ instead of $\chi$, $t^{\prime}$ denotes a time instant, the symbol $%
\left( {}\right) ,_{K_{i}}$ denotes the partial derivatives with respect to $%
K_{i}-th$ coordinate of $X,$ and is assumed the convention of sums over
repeated indices) as follows:

``Suppose that in a region $V_{0}\subset V,$
appropriate to each material body, the independent variables admit Taylor
series expansions in $X^{\prime}-X$ in $V_{0}$ [...] terminating with
gradients of order $P,Q,$ etc., 
\begin{align*}
x(X^{\prime},t^{\prime}) & =x(t^{\prime})+\left( X_{K_{1}}^{\prime
}-X_{K_{1}}\right) x,_{K_{1}}(t^{\prime}) \\
& +...+\frac{1}{P!}\left( X_{K_{1}}^{\prime}-X_{K_{1}}\right) ....\left(
X_{K_{P}}^{\prime}-X_{K_{P}}\right) x,_{K_{1}...K_{P}}(t^{\prime}),
\end{align*}
and [...]. If the response functionals are sufficiently smooth so that they
can be approximated by the functionals in the field of real functions 
\begin{align}
& x(t^{\prime}),x,_{K_{1}}(t^{\prime}),....,x,_{K_{1}...K_{P}}(t^{\prime }),
\tag{3.1.6} \\
& [...]  \notag
\end{align}
we say that the material at $X$ [...] satisfies a \textit{smooth
neighborhood hypothesis. Materials of this type, for }$P>1,Q>1$ are called 
\textit{nonsimple materials of gradient type.''}

Actually Piola is not truncating the series and keeps calculating the
integrals on the whole body $\mathcal{B}$. Although no explicit mention can
be found in the text of Piola, because of the arguments presented in remark %
\ref{BoundednessandAttenuation}, it is clear that he uses a weaker form of
the \textit{Attenuating Neighborhood Hypotheses }stated on page 34 of \cite%
{eringen2002}.

To be persuaded of this statement the reader will need to proceed to the
next section.

To conclude this section we need to remark (see Appendix C) that in very
recent times, as a karstic river, the ideas of Piola are back on the stage
of Continuum Mechanics.

The idea of an internal interaction which does not fall in the framework of
Cauchy continuum mechanics is again attracting the attention of many
researchers. Following Piola's original ideas modern "peridynamics"\footnote{%
We remark that (luckily!) the habit of inventing new names (alhough
sometimes the related concepts are not so novel) is not lost in modern
science (see \cite{russo2003} for a discussion of the importance of this
attitude in science) and that the tradition of using Greek roots for
inventing new names is still alive.} assumes that the force applied on a
material particle of a continuum actually depends on the deformation state
of a whole neighbourhood of the particle.

\subsection{An explicit calculation of the Strong Form of the Variational
Principle (\protect\ref{WORKBALANCE}).}

A more detailed discussion about the eventual novelties contained in the
formulation of peridynamics when compared with e.g. Eringen's non-local
continuum mechanics is postponed to further investigations. In this section
we limit ourselves to compute explicitly the Euler-Lagrange equation
corresponding to the Variational Principle (\ref{WORKBALANCE}). To this end
we need to treat algebraically the expression

\begin{equation}
\int_{\mathcal{B}}\left( \int_{\mathcal{B}}\Lambda(X,\bar{X},\rho)\delta
\rho^{2}\mu(\bar{X})d\bar{X}\right) \mu(X)dX  \tag{N2}  \label{NONLOCALWORK}
\end{equation}
by calculating explicitly 
\begin{equation*}
\delta\rho^{2}=\delta\left( \sum\limits_{i=1}^{3}\left( \chi_{i}(\bar {X}%
)-\chi_{i}(X)\right) \left( \chi_{i}(\bar{X})-\chi_{i}(X)\right) \right)
\end{equation*}

With simple calculations we obtain that (Einstein convention is applied from
now on)%
\begin{equation*}
\delta \rho ^{2}=\left( 2\left( \chi ^{i}(\bar{X})-\chi ^{i}(X)\right)
\left( \delta \chi _{i}(\bar{X})-\delta \chi _{i}(X)\right) \right)
\end{equation*}%
which once placed in (\ref{NONLOCALWORK}) produces 
\begin{gather*}
\int_{\mathcal{B}}\int_{\mathcal{B}}\left( 2\Lambda (X,\bar{X},\rho )\mu (%
\bar{X})\mu (X)\left( \chi ^{i}(\bar{X})-\chi ^{i}(X)\right) \right) \left(
\delta \chi _{i}(\bar{X})-\delta \chi _{i}(X)\right) d\bar{X}dX= \\
=\frac{1}{2}\left( \int_{\mathcal{B}}f^{i}(\bar{X})\delta \chi _{i}(\bar{X})d%
\bar{X}+\int_{\mathcal{B}}f^{i}(X)\delta \chi _{i}(X)dX,\right)
\end{gather*}%
where we have introduced the internal interaction force (recall that Piola,
and we agree with his considerations as presented in his num.72 on pages
150-151, assumes that $\Lambda (X,\bar{X},\rho )=\Lambda (\bar{X},X,\rho )$)
by means of the definition%
\begin{equation*}
f^{i}(\bar{X}):=\int_{\mathcal{B}}\left( 4\Lambda (X,\bar{X},\rho )\mu (\bar{%
X})\mu (X)\left( \chi ^{i}(\bar{X})-\chi ^{i}(X)\right) \right) dX
\end{equation*}%
By a standard localization argument one easily gets that (\ref{WORKBALANCE})
implies%
\begin{equation}
a^{i}(X)=b_{m}^{i}(X)+f^{i}(X)  \tag{N3}  \label{Peridynamics}
\end{equation}%
which (see also Appendix C) is exactly the starting point of modern
"peridynamics".

Many non-local continuum theories were formulated since the first
formulation by Piola: we cite here for instance 
\cite{eringen1999}, \cite{eringen2002}, \cite{eringenedelen1972}, \cite{soubestreboutin2012}.
Remarkable also are the following more modern papers \cite{dellisolaetal1995}%
, \cite{dellisolaetal2000}, \cite{dellisolaseppecher1995}, \cite%
{demmiesilling2007}, \cite{duetal2013}, \cite{emmrichetal2013}, \cite%
{lehoucqsilling2008}, \cite{selesonetal2013}, \cite{sillingetal2007}, \cite{steinmann2008}, \cite%
{steinmannetal2007}, \cite{sunyksteinmann2003}.

The non-local interaction described by the integral operators introduced in
the present subsections are not to be considered exclusively as interactions
of a mechanical nature: indeed recently a model of biologically driven
tissue growth has been introduced (see e.g. \cite{andreausetal2014}, \cite%
{madeoetal2012}, \cite{madeoetal2011}) where such a non-local operator is
conceived to model the biological stimulus to growth.

\subsection{Piola's higher gradient continua}

The state of deformation of a continuum in the neighbourhood of one of its
material points can be approximated by means of the Green deformation
measure and of all its derivatives with respect to Lagrangian referential
coordinates. Piola never considers the particular case of linearized
deformation measures (which is physically rather unnatural): his spirit has
been recovered in many modern works, among which we cite \cite%
{sillinglehoucq2008}, \cite{steigmann2002}.

Indeed in Capo VI, on page 152, Piola develops in Taylor series $\delta
\rho^{2}$ (also by using his regularity assumptions about the function $%
\Lambda(X,\bar{X},\rho)$ and the definition (11)) and replaces the obtained
development in (\ref{INTERNAL WORKATX}).

In a more modern notation (see Appendix A for the word by word translation)
starting from 
\begin{equation*}
\chi _{i}(\bar{X})-\chi _{i}(X)=\sum\limits_{N=1}^{\infty }\frac{1}{N!}%
\left( \frac{\partial ^{N}\chi _{i}(X)}{\partial X_{i_{1}}....\partial
X_{i_{N}}}(\bar{X}_{i_{1}}-X_{i_{1}})....(\bar{X}_{i_{N}}-X_{i_{N}})\right)
\end{equation*}%
Piola gets an expression for the Taylor expansion with respect to the
variable $\bar{X}$ of center $X$ for the function,%
\begin{equation*}
\rho ^{2}(\bar{X},X)=\left( \chi ^{i}(\bar{X})-\chi ^{i}(X)\right) \left(
\chi _{i}(\bar{X})-\chi _{i}(X)\right)
\end{equation*}%
He estimates and explicitly writes first, second and third derivatives of $%
\rho ^{2}$ with respect to the variable $\bar{X}$. This is what we will do
in the sequel, repeating his algebraic procedure with the only difference
consisting in the use of Levi-Civita tensor notation.

We start with the first derivative

\begin{equation}
\frac{1}{2}\frac{\partial\rho^{2}(\bar{X},X)}{\partial\bar{X}_{\alpha}}%
=\left( \chi^{i}(\bar{X})-\chi^{i}(X)\right) \frac{\partial\chi_{i}(\bar {X})%
}{\partial\bar{X}_{\alpha}}  \tag{N4}  \label{FirstDerivative}
\end{equation}
We remark that when $\bar{X}=X$ \ this derivative vanishes. Therefore the
first tem of Taylor series for $\rho^{2}$ vanishes. We now proceed by
calculating the second and third order derivatives :

\begin{align*}
\frac{1}{2}\frac{\partial^{2}\rho^{2}(\bar{X},X)}{\partial\bar{X}_{\alpha
}\partial\bar{X}_{\beta}} & =\frac{\partial\chi^{i}(\bar{X})}{\partial \bar{X%
}_{\beta}}\frac{\partial\chi_{i}(\bar{X})}{\partial\bar{X}_{\alpha}}+\left(
\chi^{i}(\bar{X})-\chi^{i}(X)\right) \frac{\partial^{2}\chi _{i}(\bar{X})}{%
\partial\bar{X}_{\alpha}\partial\bar{X}_{\beta}}= \\
& =:C_{\alpha\beta}(\bar{X})+\left( \chi^{i}(\bar{X})-\chi^{i}(X)\right) 
\frac{\partial^{2}\chi_{i}(\bar{X})}{\partial\bar{X}_{\alpha}\partial\bar {X}%
_{\beta}};
\end{align*}

\begin{equation}
\frac{1}{2}\frac{\partial^{3}\rho^{2}(\bar{X},X)}{\partial\bar{X}_{\alpha
}\partial\bar{X}_{\beta}\partial\bar{X}_{\gamma}}=\frac{\partial
C_{\alpha\beta}(\bar{X})}{\partial\bar{X}_{\gamma}}+\frac{\partial\chi _{i}(%
\bar{X})}{\partial\bar{X}_{\gamma}}\frac{\partial^{2}\chi^{i}(\bar{X})}{%
\partial\bar{X}_{\alpha}\partial\bar{X}_{\beta}}+\left( \chi^{i}(\bar {X}%
)-\chi^{i}(X)\right) \frac{\partial^{3}\chi_{i}(\bar{X})}{\partial\bar {X}%
_{\alpha}\partial\bar{X}_{\beta}\partial\bar{X}_{\gamma}}  \tag{N5}
\label{THIRDDERIVATIVE}
\end{equation}

The quantities of this last equation are exactly those described in \cite%
{piola1845-6} on page 157 concerning the quantities appearing in formulas
(14) on page 153.

We now introduce the result (formula (\ref{objectivityC})) found in Appendix
D (in order to remain closer to Piola's presentation we refrain here from
using Levi-Civita alternating symbol) 
\begin{equation*}
F_{i\gamma}\frac{\partial^{2}\chi^{i}}{\partial X^{\alpha}\partial X^{\beta}}%
=\frac{1}{2}\left( \frac{\partial C_{\alpha\gamma}}{\partial X^{\beta}}+%
\frac{\partial C_{\beta\gamma}}{\partial X^{\alpha}}-\frac{\partial
C_{\beta\alpha}}{\partial X^{\gamma}}\right)
\end{equation*}
so that by replacing in (\ref{THIRDDERIVATIVE}) we get

\begin{equation}
\frac{1}{2}\frac{\partial^{3}\rho^{2}(\bar{X},X)}{\partial\bar{X}_{\alpha
}\partial\bar{X}_{\beta}\partial\bar{X}_{\gamma}}=\frac{1}{2}\left( \frac{%
\partial C_{\alpha\gamma}}{\partial X^{\beta}}+\frac{\partial C_{\beta\gamma}%
}{\partial X^{\alpha}}+\frac{\partial C_{\beta\alpha}}{\partial X^{\gamma}}%
\right) +\left( \chi^{i}(\bar{X})-\chi^{i}(X)\right) \frac{%
\partial^{3}\chi_{i}(\bar{X})}{\partial\bar{X}_{\alpha}\partial\bar {X}%
_{\beta}\partial\bar{X}_{\gamma}}  \tag{N6}  \label{THIRDDERIVATIVELAST}
\end{equation}
so that when $\bar{X}=X$ \ we get that the third order derivatives of $%
\rho^{2}$ can be expressed in terms of the first derivatives of $%
C_{\gamma\beta}.$

Now we go back to read in Capo VI n.73 page 152-153:

`` 73. What remains to be done in order to deduce
useful consequences from the equation (12) is simply a calculation process.
Once recalled the equation (8), it is seen, transforming into series the
functions in the brackets, so that one has

\begin{align*}
\rho^{2} & =\left( f\frac{\text{ }dx}{da}+g\frac{dx}{db}+k\frac{dx}{dc}+%
\frac{f\text{ }^{2}}{2}\frac{d^{2}x}{da^{2}}+ec.\right) ^{2} \\
& +\left( f\frac{\text{ }dy}{da}+g\frac{dy}{db}+k\frac{dy}{dc}+\frac{f\text{ 
}^{2}}{2}\frac{d^{2}y}{da^{2}}+ec.\right) ^{2} \\
& +\left( f\frac{\text{ }dz}{da}+g\frac{dz}{db}+k\frac{dz}{dc}+\frac{f\text{ 
}^{2}}{2}\frac{d^{2}z}{da^{2}}+ec.\right) ^{2};
\end{align*}
and by calculating the squares and gathering the terms which have equal
coefficients:

\begin{align}
\rho ^{2}& =f\text{ }^{2}t_{1}+g\text{ }^{2}t_{2}+k\text{ }%
^{2}t_{3}+2fgt_{4}+2fkt_{5}+2gkt_{6}  \notag \\
& +f\text{ }^{3}T_{1}+2f\text{ }^{2}gT_{2}+2f\text{ }^{2}kT_{3}+f\text{ }%
g^{2}T_{4}+ec.  \tag{13}  \label{CAPOVI-13ENG2}
\end{align}%
in which expression the quantities $t_{1},t_{2},t_{3},t_{4},t_{5},t_{6}$
represent the six trimonials which are alreay familiar to us, as we have
adopted such denominations since the equations (6) in the num$%
{{}^\circ}%
$.34.; and the quantities $T_{1},T_{2},T_{3},T_{4},ec.$ where the index goes
to infinity, represent trinomials of the same nature in which derivatives of
higher and higher order appear. ''

Then the presentation of Piola continues with the study of the algebraic
structure of the trinomial constituting the quantities $T_{1},T_{2},T_{3},$
as shown by the formulas appearing in Capo VI, n.73 on pages 153-160. The
reader will painfully recognize that these huge component-wise formulas
actually have the same structure which becomes easily evident in formula \ref%
{THIRDDERIVATIVELAST} and in all formulas deduced, with Levi-Civita Tensor
Calculus, in Appendices D and E.

What Piola manages to recognize (also with a courageous conjecture, see
Appendices D and E) is that in the expression of Virtual Work all the
quantities which undergo infinitesimal variation (which are naturally to be
chosen as "measures of deformation" ) are indeed either components of the
deformation measure $C$ or components of one of its gradients.

Indeed in the num.74 page 156 one reads:

``74. A new proposition, to which the reader
should pay much attention, is that all the trinomials $%
T_{1},T_{2},T_{3},etc. $ where the index goes to infinity , which appear in
the previous equation (17), can be expresses by means of the only first six $%
t_{1},t_{2},t_{3},$ $t_{4},t_{5},t_{6,}$ and of their derivatives with
respect to the variables $a,b,c$ of all orders. I started to suspect this
analytical truth because of the necessary correspondence which must hold
between the results which are obtained with the way considered in this Capo
and those results obtained with the way considered in the Capos III and IV. 
''

This statement is true and its importance is perfectly clear to Piola: for a
discussion of the mathematical rigor of his proof the reader is referred to
the Appendix E.

In order to transform the integral expression (\ref{INTERNAL WORKATX}) 
\begin{equation*}
\left( \int_{\mathcal{B}}\Lambda(X,\bar{X},\rho)\delta\rho^{2}(X,\bar{X})\mu(%
\bar{X})d\bar{X}\right)
\end{equation*}
Piola remarks that (pages 155-156)

``When using the equation (13) to deduce the value
of the variation $\delta \rho ^{2}$ , it is clear that the characteristic $%
\delta $ will need to be applied only to the trinomials we have discussed up
to now, so that we will have:

\begin{align}
\delta \rho ^{2}& =f\text{ }^{2}\delta t_{1}+g\text{ }^{2}\delta t_{2}+k%
\text{ }^{2}\delta t_{3}+2fg\text{ }\delta t_{4}+2fk\text{ }\delta t_{5}+2gk%
\text{ }\delta t_{6}  \notag \\
& +f\text{ }^{3}\delta T_{1}+2f\text{ }^{2}g\delta T_{2}+2f\text{ }%
^{2}k\delta T_{3}+f\text{ }g^{2}\delta T_{4}+ec.  \tag{16}
\label{CAPOVI-16ENG2}
\end{align}%
Indeed the coefficients $f$ $^{2},g$ $^{2},k$ $^{2},2fg$, etc. are always of
the same form as the functions giving the variables $x,y,z$ in terms of the
variables $a,b,c,$ and therefore cannot be affected by that operation whose
aim is simply to change the form of these functions. Vice versa, by
multiplying the previous equation (16) times $\Lambda $ and then integrating
with respect to the variables $f,g,k$ \ in order to deduce from such
calculation the value to be given to the forth term under the triple
integral of the equation (12), such an operation is affecting only the
quantities $\Lambda f$ $^{2},\Lambda g$ $^{2}$, etc. and the variations $%
\delta t_{1},\delta t_{2},\delta t_{3}....\delta T_{1},\delta T_{2},ec.$
cannot be affected by it, as the trinomials $%
t_{1},t_{2},t_{3}....T_{1},T_{2},ec.$ (one has to consider carefully which
is their origin) do not contain the variables $f,g,k$ : therefore such
variations result to be constant factors, times which are to be multiplied
the integrals to be calculated in the subsequent terms of the series. 
''

Using a modern notation we have that 
\begin{equation*}
\rho ^{2}(\bar{X},X)=\sum\limits_{N=1}^{\infty }\frac{1}{N!}\left. \frac{%
\partial ^{N}\rho ^{2}(\bar{X},X)}{\partial \bar{X}_{i_{1}}....\partial \bar{%
X}_{i_{N}}}\right\vert _{X=\bar{X}}(\bar{X}_{i_{1}}-X_{i_{1}})....(\bar{X}%
_{i_{N}}-X_{i_{N}})
\end{equation*}%
and therefore that 
\begin{equation*}
\delta \rho ^{2}(\bar{X},X)=\sum\limits_{N=1}^{\infty }\frac{1}{N!}\left(
\delta \left. \frac{\partial ^{N}\rho ^{2}(\bar{X},X)}{\partial \bar{X}%
_{i_{1}}....\partial \bar{X}_{i_{N}}}\right\vert _{X=\bar{X}}\right) (\bar{X}%
_{i_{1}}-X_{i_{1}})....(\bar{X}_{i_{N}}-X_{i_{N}}).
\end{equation*}%
As a consequence 
\begin{gather*}
\int_{\mathcal{B}}\Lambda (X,\bar{X},\rho )\delta \rho ^{2}(\bar{X},X)\mu (%
\bar{X})d\bar{X}= \\
=\sum\limits_{N=1}^{\infty }\frac{1}{N!}\left( \delta \left. \frac{\partial
^{N}\rho ^{2}(\bar{X},X)}{\partial \bar{X}^{i_{1}}....\partial \bar{X}%
^{i_{N}}}\right\vert _{X=\bar{X}}\right) \left( \int_{\mathcal{B}}\Lambda (X,%
\bar{X},\rho )\left( (\bar{X}^{i_{1}}-X^{i_{1}})....(\bar{X}%
^{i_{N}}-X^{i_{N}})\right) \mu (\bar{X})d\bar{X}\right)
\end{gather*}

If we introduce the tensors 
\begin{equation*}
T_{.}^{i_{1}...i_{N}}(X):=\left( \int_{\mathcal{B}}\Lambda (X,\bar{X},\rho
)\left( (\bar{X}^{i_{1}}-X^{i_{1}})....(\bar{X}^{i_{N}}-X^{i_{N}})\right)
\mu (\bar{X})d\bar{X}\right)
\end{equation*}%
we get, also by recalling formula (\ref{HYPREC1}) from Appendix E, 
\begin{equation*}
\int_{\mathcal{B}}\Lambda (X,\bar{X},\rho )\delta \rho ^{2}(\bar{X},X)\mu (%
\bar{X})d\bar{X}=\sum\limits_{N=1}^{\infty }\frac{1}{N!}\left( \delta
L_{\alpha _{1}....\alpha _{n}}\left( C(X),..,\nabla ^{n-2}C(X)\right)
\right) T_{.}^{i_{1}...i_{N}}(X)
\end{equation*}

Piola states that

``After these considerations it is manifest the
truth of the equation:

\begin{equation}
\int df\int dg\int dk\cdot \Lambda \delta \rho ^{2}=  \tag{17}
\label{CAPOVI-17ENG2}
\end{equation}

\begin{align*}
& \left( 1\right) \text{ }\delta t_{1}+\left( 2\right) \text{ }\delta
t_{2}+\left( 3\right) \text{ }\delta t_{3}+\left( 4\right) \text{ }\delta
t_{4}+\left( 5\right) \text{ }\delta t_{5}+\left( 6\right) \text{ }\delta
t_{6} \\
& +\left( 7\right) \text{ }\delta T_{1}+\left( 8\right) \text{ }\delta
T_{2}+\left( 9\right) \text{ }\delta T_{3}+\left( 10\right) \text{ }\delta
T_{4}+ec.
\end{align*}
where the coefficients (1), (2), etc. indicated by means of numbers in
between brackets, must be regarded to be each a function of the variables $%
a,b,c$ as obtained after having performed the said definite integrals. 
''

In order to establish the correct identification between Piola's notation
and the more modern notation which we have introduced the reader may simply
consider the following table ($i=1,2,....n,....)$%
\begin{equation*}
\begin{array}{cc}
T_{.}^{i_{1}...i_{N}}\leftrightarrows(1),(2),ec. & \delta L_{\alpha
_{1}....\alpha_{n}}\left( C,..,\nabla^{n-2}C\right) \leftrightarrows\delta
T_{i}%
\end{array}
.
\end{equation*}

After having accepted Piola's assumptions the identity (\ref{WORKBALANCE})
becomes 
\begin{multline}
\int_{\mathcal{B}}\left( \left( b_{m}(X)-a(X)\right) \delta \chi
(X)+\sum\limits_{N=1}^{\infty }\frac{1}{N!}\left( \delta L_{\alpha
_{1}....\alpha _{n}}\left( C(X),..,\nabla ^{n-2}C(X)\right) \right)
T_{.}^{i_{1}...i_{N}}(X)\right) \mu (X)dX  \label{N-thGradient} \\
+\delta W(\partial \mathcal{B})=0  \notag
\end{multline}%
By a simple re-arrangement and by introducing a suitable notation the last
formula becomes 
\begin{equation}
\int_{\mathcal{B}}\left( \left( b_{m}(X)-a(X)\right) \delta \chi
(X)+\sum\limits_{N=1}^{\infty }\left\langle \nabla ^{N}\delta
C(X)|S_{.}(X)\right\rangle \right) \mu (X)dX+\delta W(\partial \mathcal{B})=0
\tag{12tris}  \label{N-thGradientFinal}
\end{equation}%
where $S$ is a $N-th$ order contravariant totally symmetric tensor\footnote{%
The constitutive equations for such tensors must verify the condition of
frame invariance. When these tensors are defined in terms of a deformation
energy (that is when the Principle of Virtual Work is obtained as the first
variation of a Least Action Principle) the objectivity becomes a restriction
on such an energy. The generalization of the results in Steigmann (2003) to
the N-the gradient continua still needs to be found.} and the symbol $%
\left\langle 
\begin{array}{c}
\end{array}%
|%
\begin{array}{c}
\end{array}%
\right\rangle $ denotes the total saturation (inner product) of a pair of
totally symmetric contravariant and covariant tensors.

Indeed on pages 159-160 of \cite{piola1845-6} we read

``75. Once the proposition of the previous num.
has been admitted, it is manifest that the equation (17) can assume the
following other form 
\begin{equation}
\int df\int dg\int dk\cdot \Lambda \delta \rho ^{2}=  \tag{18}
\label{CAPOVI-18ENG2}
\end{equation}

\begin{align*}
& \left( \alpha \right) \text{ }\delta t_{1}+\left( \beta \right) \text{ }%
\delta t_{2}+\left( \gamma \right) \text{ }\delta t_{3}+....+\left( \epsilon
\right) \text{ }\frac{\delta dt_{_{1}}}{da}+\left( \zeta \right) \text{ }%
\frac{\delta dt_{_{1}}}{db}+\left( \eta \right) \text{ }\frac{\delta
dt_{_{1}}}{dc} \\
& +\left( \vartheta \right) \text{ }\frac{\delta dt_{_{2}}}{da}+....+\left(
\lambda \right) \text{ }\frac{\delta d^{2}t_{_{1}}}{da^{2}}+\left( \mu
\right) \text{ }\frac{\delta d^{2}t_{_{1}}}{dadb}+....+\left( \xi \right) 
\text{ }\frac{\delta d^{2}t_{_{2}}}{da^{2}}+\left( o\right) \frac{\delta
d^{2}t_{_{2}}}{dadb}+ec.
\end{align*}%
in which the coefficients $\left( \alpha \right) ,\left( \beta \right) $ $%
....\left( \epsilon \right) ....\left( \lambda \right) ....ec.$ \ are
suitable quantities given in terms of the coefficients $\left( 1\right)
,\left( 2\right) ....\left( 7\right) ,\left( 8\right) ....$ of the equation
(17): they depend on the quantities $t_{1},t_{2}....t_{6}$, and on all order
derivatives of these trinomials with respect to the variables $a,b,c$ . Then
the variations $\delta t_{1},$ $\delta t_{2}...$.(with the index varying up
to infinity) and the variations of all their derivatives of all orders $%
\frac{\delta dt_{_{1}}}{da},\frac{\delta dt_{_{1}}}{db},ec.$ appear in the
(18) only linearly ''

\section{Weak and Strong Evolution Equations for Piola Continua}

To our knowledge a formulation of the Principle of Virtual Work for $N-th$
gradient Piola Continua equivalent to (\ref{N-thGradientFinal}) is found in
the literature only in \cite{dellisolaetal2012}, but the authors were
unaware of the previous work of Piola.

The reader is referred to the aforementioned paper for the detailed
presentation of the needed Postulation process and the subsequent procedure
of integration by parts needed for transforming the weak formulation of
evolution equations given by (\ref{N-thGradientFinal}) into a strong
formulation in which suitable bulk equations and the corresponding boundary
conditions are considered.

We shortly comment here about the relative role of Weak and Strong
formulations, framing it in a historical perspective.

Since at least the pioneering works by Lagrange the Postulation process for
Mechanical Theories was based on the Least Action Principle or on the
Principle of Virtual Work.

One can call "Variational" both these Principles as the Stationarity
Condition for Least Action requires that for all admissible variations of
motion the first variation of Action must vanish, statement which, as
already recognized by Lagrange him-self, implies a form of the Principle of
Virtual Work.

However in order to "compute" the motions relative to given initial data the
initiators of Physical Theories needed to integrate by parts the
Stationarity Condition which they had to handle.

In this way they derived some PDEs with some boundary conditions which
sometimes were solved by using analytical or semi-analytical methods.

From the mathematical point of view this procedure is applicable when the
searched solution have a stronger regularity than the one strictly needed to
formulate the basic variational principle.

It is a rather ironic circumstance that nowadays very often those
mathematicians who want to prove well-posedness theorems for PDEs (which
originally were obtained by means of an integration by part procedure) start
their reasonings by applying in the reverse direction the same integration
by parts process: indeed very often the originating variational principle of
all PDEs is forgotten. Some examples of mathematical results which exploit
in an efficient way the power of variational methods are those presented for
instance in Neff \cite{neff2007}, \cite{neff2002}, \cite{neff2006c}.

Actually even if one refuses to accept the idea of basing all physical
theories on variational principles, he is indeed obliged, in order to find
the correct mathematical frame for his models, to "prove" the validity of a
weak form appliable to his painfully formulated balance laws. In reality
(see \cite{dellisolaplacidi2011}) his model will not be acceptable until he
has been able to reformulate it in a "weak" form.

It seems that the process which occurred in mathematical geography,
described in Russo \cite{russo2013}-\cite{russo2003}, occurs very frequently
in science. While the reader is referred to the cited works for all details,
we recall here the crucial point of Russo's argument, as needed for our
considerations. Ptolemy presented in his Almagest a useful tool for
astronomical calculations: actually his \textit{Handy Tables} tabulate all
the data needed to compute the positions of the Sun, Moon and planets, the
rising and setting of the stars, and eclipses of the Sun and Moon. The main
calculation tools in Ptolemy's treatise are the deferents and epicycles,
which were introduced by Apollonius of Perga and Hipparchus of Rhodes in the
framework of astronomic theories much more advanced than the one formulated
by Ptolemy (if Russo's conjecture is true). Unfortunately Ptolomy
misunderstood the most ancient (and much deeper) theories and badly
re-organized the knowledges, observations and theories presented in the
treatise by Hipparchus (treatise which has been lost): indeed Ptolemy being
"a practical scientist" gives a too high importance to "the calculation
tools" by blurring in a list of logical incongruences the rigorous and deep
(and eliocentric!) theories formulated by Hipparchus four centuries before
him. Actually in Ptolemy's vision the calculation tools become the
fundamental ingredients of the mathematical model which he presents.

This seems to have occurred also in Continuum Mechanics: the Euler-Lagrange
equations, obtained by means of a process of integration by parts, were
originally written, starting from a variational principle, to supply a
"calculation tool" to applied scientists. They soon became (for simplifying)
the "bulk" of the theories and often the originating variational principles
were forgotten (or despised as too "mathematical"). For a period balance
equations were (with some difficulties which are discussed e.g. in \cite%
{dellisolaplacidi2011}) postulated "on physical grounds".

When the need of proving rigorous existence and uniqueness theorems met the
need of developing suitable numerical methods, and when the many failures of
the finite difference schemes became evident, the variational principles
were re-discovered \textit{starting from the balance equations.}

The variational principles represent at first the starting point of
mechanical theories and were used to get, by means of algebraic
manipulation, some tools for performing "practical calculations":\ i.e., the
associated Euler-Lagrange equations or (using another name) balance
equations. However, with a strange exchange of roles, if their basic role is
forgotten and balance equations are regarded as the basic principles from
which one has to start the formulation of the theories, then variational
principles need to be recovered as a computational tool.

One question needs to be answered: why in the modern paper \cite%
{dellisolaetal2012} a strong formulation was searched of the evolution
equation for $N-th$ gradient continua? The answer is: beacause of the need
of finding for those theories the most suitable boundary conditions !

This point is discussed also in \cite{piola1845-6} as remarked already in 
\cite{auffreyetal2014}.

\cite{piola1845-6} on pages 160-161 claims that

``Now it is a fundamental principle of the
calculus of variations (and we used it also in this Memoir in the num.$%
{{}^\circ}%
$ 36. and elsewhere) that one series as the previous one, where the
variations of some quantities and the variations of their derivatives with
respect to the fundamental variables $a,b,c$ appear linearly can be always
be transformed into one expression which containes the first quantities
without any sign of derivation, with the addition of other terms which are
exact derivatives with respect to one of the three simple independent
variables. As a consequence of such a principle, the expression which
follows to the equation (18) can be given

\begin{equation}
\int df\int dg\int dk\cdot \Lambda \delta \rho ^{2}=  \tag{19}
\label{CAPOVI-19ENG2}
\end{equation}

\begin{align*}
& A\delta t_{1}+B\text{ }\delta t_{2}+C\text{ }\delta t_{3}+D\delta
t_{4}+E\delta t_{5}+F\text{ }\delta t_{6} \\
& \text{ \ \ \ \ \ \ \ \ \ \ \ \ \ \ \ }+\frac{d\Delta }{da}+\frac{d\Theta }{%
db}+\frac{d\Upsilon }{dc}.
\end{align*}%
The values of the six coefficients $A,B,C,D,E,F$ are series constructed with
the coefficients $\left( \alpha \right) ,\left( \beta \right) ,\left( \gamma
\right) ....\left( \epsilon \right) ,\left( \zeta \right) ....\left( \lambda
\right) ,ec.$ of the equation (18) which appear linearly, with alternating
signs and affected by derivations of higher and higher order when we move
ahead in the terms of said series: the quantities $\Delta ,\Theta ,\Upsilon $
\ are series of the same form of the terms which are transformed, in which
the coefficients of the variations have a composition similar to the one
which we have described for the six coefficients $A,B,C,D,E,F$ .

Once -instead of the quantity under the integral sign in the left hand side
of the equation (12)- one introduces the quantity which is on the right hand
side of the equation (19), it is clear to everybody that an integration is
possible for each of the last three addends appearing in it and that as a
consequence these terms only give quantities which supply boundary
conditions. What remains under the triple integral is the only sestinomial
which is absolutely similar to the sestinomial already used in the equation
(10) num.$%
{{}^\circ}%
$ 35. for rigid systems. Therefore after having remarked the aforementioned
similarity the analytical procedure to be used here will result perfectly
equal to the one used in the num.$%
{{}^\circ}%
$ 35, procedure which led to the equations (26), (29) in the num.$%
{{}^\circ}%
$ 38 and it will become possible the demonstration of the extension of the
said equations to every kind of bodies which do not respect the constraint
of rigidity, as it was mentioned at the end of the num.$%
{{}^\circ}%
$ 38. It will also be visible the coincidence of the obtained results with
those which are expressed in the equations (23) of the num.$%
{{}^\circ}%
$50. which hold for every kind of systems and which were shown in the Capo
IV by means of those intermediate coordinates $p,q,r$, whose consideration,
when using the approach used in this Capo, will not be needed.''

The novel content in \cite{dellisolaetal2012} consists in the determination
of

\begin{itemize}
\item the exact structure of the tensorial quantity whose components are
called $A,B,C,D,E,F$ by \cite{piola1845-6}

\item the exact structure of the boundary conditions resulting when applying
Gauss' theorem to the divergence field called by \cite{piola1845-6} 
\begin{equation*}
\frac{d\Delta }{da}+\frac{d\Theta }{db}+\frac{d\Upsilon }{dc}
\end{equation*}%
on a suitable class of contact surfaces.
\end{itemize}

The considerations sketched about the history of celestial mechanics should
persuade the reader that it is not too unlikely that some ideas by Piola
needed 167 years for being further developed (even if the fact that the
authors did not manage to find any intermediate reference does not mean that
such a reference does not exist, maybe in a language even less
understandable than Italian).

Earlier papers (nowadays considered fundamental) by Mindlin \cite%
{mindlin1964}, \cite{mindlin1965}, \cite{mindlinetal1968}, \cite{sedov1972}, 
\cite{sedov1968} had developed a more complete study of Piola Continua, at
least up to those whose deformation energy depends on the Third Gradient,
completely characterizing the nature of contact actions in these cases, or
for continua having a kinematics richer than that considered by Piola,
including microdeformations and micro-rotations.

Many important applications can be conceived for higher gradient materials,
as for instance those involving the phenomena described for instance in \cite%
{alibertetal2003}, \cite{dellisolaetal2000}, \cite{dellisolaetal2009}, \cite%
{gatignolseppecher1986}, \cite{madeogavrilyuk2010}, \cite{madeoetal2012}, 
\cite{rosietal2013}, \cite{seppecher1996a}, \cite{seppecher1993}, \cite%
{seppecher1987}, \cite{seppecher1989}, \cite{seppecher1996b}, \cite%
{seppecher2001}, \cite{seppecher2002}, \cite{seppecher1988}, \cite{yangetal2011}, \cite%
{yangmisra2010}, \cite{yangmisra2012}, \cite%
{yeremeyevetal2007}.

\section{One- and Two-dimensional Continua and Micro-Macro identification
procedure as introduced by Piola (1845-6)}

On page 19 Piola justifies the introduction of one-dimensional or
two-dimensional bodies as follows

``11. Sometimes mathematicians are used to
consider the matter configured not in a volume with three dimensions but
[configured] in a line or in a surface: in these cases we have the so called
linear or surface systems. Indeed [these systems] are nothing else than
abstractions and it is just for this reason that the Geometer should pay the
major attention to three dimensional systems. Nevertheless, it is useful to
consider [these systems] because the several analyses for the three kind of
systems provide feedbacks that make clear [such analyses], and moreover
[such analyses] are useful for physical applications, eventhough always in
an approximate way, because the bodies, rigorously speaking, being never
deprived in Nature of one or two dimensions.

Although for both linear and surface systems we need special considerations
in order to represent the distribution of the molecules, and [in order] to
form the idea of the density and the measure of the mass, yet [the idea and
the measure] are at all similar to the above referred for three dimensional
systems: thus, I will expound them shortly. ''

On page 39 num. 24 and on page 46 num. 29 of \cite{piola1845-6} is studied
the structure of the Principle of Virtual Work in the case in which one or
two dimensions of the considered body can be neglected in the description of
its motion.

Piola uses these parts to prepare the reader for the micro-macro
identification process for three-dimensional bodies which he will study
later in full detail.

This identification process

\begin{itemize}
\item starts from a discrete system of material particles which are placed
in a reference configuration at the nodes of a suitably introduced mesh,

\item proceeds with the introduction of a suitable placement field $\chi$
having all the needed regularity properties

\item assumes that the values of $\chi$ at the aforementioned nodes can be
considered an approximation of the displacements of the discrete system of
material particles

\item and is based on the identification of Virtual Work expressions in the
discrete and continuous models.
\end{itemize}

While the detailed description of aforementioned identification (see \cite%
{alibertetal2003}, \cite{ataisteigmann1997}, \cite{haseganusteigmann1996})
process is postponed to further studies, we want here to remark that
non-local and higher gradient theories for beams and shells are already
implicitly formulated in \cite{piola1845-6}, although the main subject there
is the study of three-dimensional bodies.

The authors have found interesting connections in this context with many of
the subsequent works and the most suggestive are those concerning the theory
of shells and plates; namely, \cite{eremeyevetal2003}, \cite%
{eremeyevlebedev2013}, \cite{eremeyevpietraszkiewicz2004}, \cite%
{eremeyevpietraszkiewicz2011}, \cite{madeoetal2013}, \cite%
{pietraszkiewiczetal2007}, where interesting phenomena involving phase
transition are considered, or the papers by Neff \cite{neff2006a}, \cite%
{neff2004}, \cite{neff2007}, \cite{neff2006b}, \cite{neff2006c}, \cite%
{neffchelminski2007}.

Moreover the methods started by Piola are used also when describing
bidimensional surfaces carrying material properties as for instance in \cite%
{lebedevetal2010}, \cite{mcbrideetal2011}, \cite{mcbrideetal2012}, \cite%
{placidietal2013}, \cite{steebdiebels2004}, \cite{steigmann2009}, \cite%
{steigmannogden1999}, \cite{steigmannogden1997}.

Also interesting analogies for what concerns the connections between
discrete and continuous models can be found with papers dealing with
one-dimensional continua and their stability as for instance \cite%
{luongodiegidio2005} and \cite{luongoromeo2006}, where are studied the
dynamics of beams or chains of beams, \cite{luongoetal2009}, where the
non-linear equations for inextensible cables deduced by Piola are applied to
very interesting special motions, \cite{steigmann1992} where the case of
prestressed networks is considered, \cite{steigmann1996}\hspace{0in}and \cite%
{steigmannfaulkner1993}, where the spirit of Piola's contribution is adapted
to the context of spatial rods and the nonlinear theory for spatial
lattices. Concerning the micro-macro identification procedure in the recent
literature one can find many continuators of Piola's works. Notable are the
works \cite{boutinetal2003}, \cite{boutinhans2003}, \cite{chesnaisetal2012}
in which Piola continua are obtained by means of homogenization procedures
starting from lattice beam microstructures. It is possible to cite also some
studies which consider visco-elastic continuum theories with damage (see 
\cite{contrafattocuomo2006}, \cite{contrafattocuomo2005}, \cite{contrafattocuomo2002}, \cite{cuomocontrafatto2000}, \cite{cuomoventura1998}) for microscopically granular or discrete systems as for instance \cite%
{misrachang1993}, \cite{misraching2013}, \cite{misrasingh2013}, \cite%
{misrayang2010}, \cite{rinaldilai2007}, \cite{rinaldietal2008} or other
studies of phenomena involving multiscale coupling (see e.g. \cite%
{nadleretal2006}).

\section{A Conclusion: Piola as precursor of the Italian School of
Differential Geometry}

The most important contribution of Gabrio Piola to mechanical sciences is
the universally recognized Piola transformation, which allows for the
transformation of some equations in a conservative form from Lagrangian to
Eulerian description. The differential geometric content of this
contribution does not need to be discussed, as it has been treated in many
works and textbooks: we simply refer to \cite{epstein2010} and to the
references there cited for a detailed discussion of this point and more
considerations about the relationship between continuum mechanics and
differential geometry (see also \cite{epsteinsegev1980}).

In the present paper we have shown that there are other major contributions
to mechanics by Gabrio Piola which have been underestimated: we also have
tried a first analysis of the reasons for which this circumstance occurred.

In this concluding section we want to remark that also those results by
Piola which we have described in the present paper have a strong connection
with differential geometry (in this context see also \cite{segev1986}, \cite%
{segev2000}).

The readers is referred to the discussion about "historical method" which
was developped in the Introduction: knowledge of the basic ideas of
differential geometry is required to follow the considerations which we
present here. The criticism usually given to the kind of reconstructions
which we want to present is usually based on the following statement: the
historian wanted to read something which could not be written in such an
early stage of knowledge.

We dismiss a priori this criticism on the basis of the following statements

\begin{itemize}
\item The inaugural lecture by Riemann dates to 1854 therefore Piola's
results are surely antecedent but very close in time.

\item Riemann is considered one of the founders of Riemannian geometry even
if he did not write any formula using the indicial notation developed by
Ricci and Levi-Civita

\item Riemannian tensor is named after Riemann even if there is no formal
definition of the concept of tensor in Riemann's works.
\end{itemize}

In his inaugural lecture Riemann discusses one of his main contributions to
geometry: i.e. the condition for which a Riemannian manifold is flat. This
study (indirectly influenced by Gauss) started a flow of investigations in
which the Italian School has played a dominant role. We recall here e.g.
Ricci's Lemma and Identities, the concept of Levi-Civita parallel transport
and the Levi-Civita Theorem about parallel transports compatible with a
Riemannian structure. Also referring to the Appendix F for substantiating
our statement we claim that it was indeed Continuum Mechanics which
originated Differential Geometry and that the Italian School in differential
geometry may have been originated in the works of Piola. Indeed in Appendix
D we have proven that Piola has obtained (component-wise, exactly in the
same form in which Riemann obtained all his results) the equation (\ref%
{RiemannIdentity}),

\begin{equation*}
F_{i\gamma }\frac{\partial ^{2}\chi ^{i}}{\partial X^{\alpha }\partial
X^{\beta }}=\frac{1}{2}\left( \frac{\partial C_{\alpha \gamma }}{\partial
X^{\beta }}+\frac{\partial C_{\beta \gamma }}{\partial X^{\alpha }}-\frac{%
\partial C_{\beta \alpha }}{\partial X^{\gamma }}\right) .
\end{equation*}%
This equation is equivalent (see \cite{spivak1979} vol.II page 184 ) to the
Riemannian condition of flatness.

\section{Appendices}

\subsection{Appendix A. Verbatim translation of excerptions from Piola
(1845-6).}

In the following subsubsections, one page for each subsubsection, we
translate from Italian to English the page indicated in the respective title.

\subsubsection{CAPO VI \textit{On the motion of a generic [deformable] body
following the ideas of the modern scientists about the molecular actions,
page 146}}

At the beginning of the Capo IV it was said that there are two ways for
taking into account -in the general equation of the motion of a generic
body- the effect of the constraints established among its molecules by
internal forces A [first ] way which was introduced consisted in expressing
such effects by means of equations of condition, and therefore by means of
the third part of the most general equation (1) in the num$%
{{}^\circ}%
$.16. : this was the way which we used in the preceding two Chapters. A
second way consisted in considering -following the ideas of modern
scientists- the molecular actions by making use of the second part of the
aforementioned equation (1), where are to be included all the terms
introduced by internal active forces: about this second way I will discuss
now. This effort will be performed also because we will see that the two
different ways actually lead to the same results at least for that part of
the solution which is the most important and fundamental (and this agreement
is really very reassuring): on the other hand it has to be remarked that the
two said ways are completing each other, and one sheds light on the other so
that what was complicated and difficult in one way becomes easy in the other
one.

\subsubsection{Page 147}

71. Recalling what was said in the numbers 31, 32 to show how, in the case of
the systems having three dimensions, the first part of the general equation
(1) num$%
{{}^\circ}%
$.16., due to external forces, is modified to be represented as follows:

\begin{equation}
\int da\int db\int dc\cdot\left\{ \left( X-\frac{d^{2}x}{dt^{2}}\right)
\delta x+ec.\right\} ;  \tag{1}  \label{CAPOVI-1ENG}
\end{equation}
we see now how it has to be modified the second part $Sm_{i}m_{j}$ $%
K\delta\rho$ , which is that one we want to consider now, while at the same
time in the third part we equate to zero all terms expressing actions
applied to all the mass [of considered body] and only retain those terms
related to forces concentrated on surfaces, lines and points.

This second part, once assuming that for each pair of molecules there is
acting always an internal force $K$, when the number of points is equal to $%
n $, when expressed explicitly can be represented as follows:

\begin{align}
&
m_{_{1}}m_{_{2}}K_{_{1,2}}\delta\rho_{_{1,2}}+m_{_{1}}m_{_{3}}K_{_{1,3}}%
\delta\rho_{_{1,3}}+.....+m_{_{1}}m_{n}K_{_{1},n}\delta\rho_{_{1},n}  \notag
\\
& \text{ \ \ \ \ \ \ \ \ \ \ \ \ \ \ \ \ \ \ }+m_{_{2}}m_{_{3}}K_{_{2,3}}%
\delta\rho_{_{2,3}}+.....+m_{_{2}}m_{n}K_{_{2,}n}\delta\rho_{_{2,}n}  \tag{2}
\label{CAPOVI-2ENG} \\
& \text{ \ \ \ \ \ \ \ \ \ \ \ \ \ \ \ \ \ \ \ \ \ \ \ \ \ \ \ \ \ \ \ \ \ \
\ \ \ \ \ \ \ \ \ \ \ \ \ \ \ \ \ \ }\vdots\text{ \ \ \ \ \ \ \ \ \ \ \ \ }%
\vdots  \notag \\
& \text{ \ \ \ \ \ \ \ \ \ \ \ \ \ \ \ \ \ \ \ \ \ \ \ \ \ \ \ \ \ \ \ \ \ \
\ \ \ \ \ \ \ \ \ \ \ \ \ \ \ \ \ \ \ \ \ \ \ \ \ }+m_{n-1}m_{_{n}}K_{n-1,n}%
\delta\rho_{n-1,n}  \notag
\end{align}
being in general:

\begin{equation}
\rho _{i,j}=\sqrt{\left( x_{j}-x_{i}\right) ^{2}+\left( y_{j}-y_{i}\right)
^{2}+\left( z_{j}-z_{i}\right) ^{2}.}  \tag{3}  \label{CAPOVI-3ENG}
\end{equation}%
It can be however seen that the subsequent horizontal lines appearing in it
, which one after another have a lacking term with respect to the previous
line, can be rewritten in such a way that all have exactly $n$ terms, by
writing, at the place of the equation (2) the quantity

\begin{align}
\frac{1}{2} & m_{_{1}}m_{_{1}}K_{_{1,1}}\delta\rho_{_{1,1}}+\frac{1}{2}%
m_{_{1}}m_{_{2}}K_{_{1,2}}\delta\rho_{_{1,2}}+.....+\frac{1}{2}%
m_{_{1}}m_{n}K_{_{1},n}\delta\rho_{_{1},n}  \notag \\
& +\frac{1}{2}m_{_{2}}m_{_{1}}K_{_{2,1}}\delta\rho_{_{2,1}}+\frac{1}{2}%
m_{_{2}}m_{_{2}}K_{_{2,2}}\delta\rho_{_{2,2}}+.....+\frac{1}{2}%
m_{_{2}}m_{n}K_{_{2},n}\delta\rho_{_{2},n}  \notag \\
& \text{\ \ \ \ \ \ \ \ \ \ \ \ \ \ \ \ \ \ \ \ \ \ \ }\vdots\text{ \ \ \ \
\ \ \ \ \ \ \ \ \ \ \ \ \ }\vdots  \tag{4}  \label{CAPOVI-4ENG} \\
& +\frac{1}{2}m_{i}m_{_{1}}K_{i,_{1}}\delta\rho_{i,_{1}}+\frac{1}{2}%
m_{i}m_{_{2}}K_{i,2}\delta\rho_{i,2}+.....+m_{i}m_{n}K_{i,n}\delta\rho _{i,n}
\notag \\
& \text{ \ \ \ \ \ \ \ \ \ \ \ \ \ \ \ \ \ \ \ \ \ \ \ }\vdots\text{ \ \ \ \
\ \ \ \ \ \ \ \ \ \ }\vdots  \notag \\
& +\frac{1}{2}m_{n}m_{_{1}}K_{n,_{1}}\delta\rho_{n,_{1}}+\frac{1}{2}%
m_{n}m_{_{2}}K_{n,2}\delta\rho_{n,2}+.....+\frac{1}{2}m_{n}m_{n}K_{n,n}%
\delta\rho_{n,n}.  \notag
\end{align}

\subsubsection{Page 148}

\noindent To recognize the equality of the two quantities (2), (4), it is
enough to observe first that in the second one the terms containing the
variations

\begin{equation*}
\delta \rho _{_{1,1}},\delta \rho _{_{2,2}}.....\delta \rho
_{i,i}.....\delta \rho _{n,n}
\end{equation*}%
are introduced only for maintaining the regularity in the progression of the
indices, and it is as if they where not present, as the radicals $\rho
_{_{1,1}},\rho _{_{2,2}}....\rho _{n,n},$ and their variations are
vanishing, as it is manifest from the generic expression (3). Secondly it
has to be observed that the remaining terms can be pair-wise added:
therefore the two terms $\frac{1}{2}m_{_{1}}m_{_{2}}K_{_{1,2}}\delta \rho
_{_{1,2}}+\frac{1}{2}m_{_{2}}m_{_{1}}K_{_{2,1}}\delta \rho _{_{2,1}}$ are
equivalent to the following one $m_{_{1}}m_{_{2}}K_{_{1,2}}\delta \rho
_{_{1,2}}$ . Indeed it is clear that, because of (3) the quantity $\rho
_{_{1,2}}$ equals the quantity $\rho _{_{2,1}}$ : and that the force $%
K_{_{1,2}}$ equals the force $K_{_{2,1}}$ , as it is implied by the
Principle of Action and Reaction and as it will become even clearer for what
we will add about the way in which the generic composition of the internal
action $K$ has to be understood. In a similar way the two terms $\frac{1}{2}%
m_{_{1}}m_{_{3}}K_{_{1,3}}\delta \rho _{_{1,3}}+\frac{1}{2}%
m_{_{3}}m_{_{1}}K_{_{3,1}}\delta \rho _{_{3,1}}$ will gather into a single
one $m_{_{1}}m_{_{3}}K_{_{1,3}}\delta \rho _{_{1,3}}$ ; and so on for the
other terms. After all previous observations it is easy to persuade oneself
that the quantity (4) is equivalent to the shorter form given by (2).

Any whatsoever of the horizontal series which compose the quantity (4) can
be reduced by means of a triple summation. To be persuaded of the truth of
this statement it is needed to recall the idea of the previously introduced
disposition of the molecules by means of the coordinates $a,b,c$ which
allows us to represent the coordinates of the generic molecule $m_{i}$ as
given by the following functions

\begin{equation}
x_{i}=x\left( a,b,c\right) ;\text{ \ \ \ }y_{i}=y\left( a,b,c\right) ;\text{
\ \ \ }z_{i}=z\left( a,b,c\right) .  \tag{5}  \label{CAPOVI-5ENG}
\end{equation}
Once we have also represented as follows:

\begin{equation}
x_{j}=x\left( a+f,b+g,c+k\right) ;\text{ \ \ \ }y_{j}=y\left(
a+f,b+g,c+k\right) ;\text{ \ \ \ }z_{j}=z\left( a+f,b+g,c+k\right)  \tag{6}
\label{CAPOVI-6ENG}
\end{equation}%
the coordinates $x_{j},y_{j},$\ $z_{j}$ of the other whatsoever molecule $%
m_{j}$ , which (if the molecule $m_{i}$ is kept as fixed) will subsequently
pass to mean all the other molecules of considered mass; and we mean that
these analytical values (6) will vary following the variation of the
molecules $m_{j}$ when in them the variables $f,g,k,$ are changed, while the
variables $a,b,c$. are kept fixed. This is done as if we were imagining in
the preceding ideal configuration that three new axes having as origin the
molecule $m_{i}$ and parallel to those relative to the variables $a,b,c,$
have been introduced,

\subsubsection{Page 149}

\noindent and as if we were calling$,$ $f,g,k$ the coordinates of a molecule
whatsoever with respect to said new axes. Now it is convenient to recall
what was said in the num$%
{{}^\circ}%
$.31., when the first part of the general equation was treated, about the
way in which one can imagine how the $n$ points of the system are
distributed following the there given positions relatively to the three
axes, which lead to give to the ensemble of $n$ terms [appearing in that
general equation] a structure of triple series: and it will not be difficult
to understand that the $(i)-th$ of the horizontal series which compose the
quantity (4) can be represented by means of the following finite triplicate
integral

\begin{equation}
\Sigma \Delta f\Sigma \Delta g\Sigma \Delta k\cdot \frac{1}{2}%
m_{i}m_{j}K\delta \rho ,  \tag{7}  \label{CAPOVI-7ENG}
\end{equation}%
where the quantity $\rho $ (equations (3),(5),(6)) has the value given by
the equation

\begin{align}
\rho ^{2}& =\left[ x\left( a+f,b+g,c+k\right) -x\left( a,b,c\right) \right]
^{2}  \notag \\
& +\left[ y\left( a+f,b+g,c+k\right) -y\left( a,b,c\right) \right] ^{2} 
\tag{8}  \label{CAPOVI-8ENG} \\
& +\left[ z\left( a+f,b+g,c+k\right) -z\left( a,b,c\right) \right] ^{2}. 
\notag
\end{align}%
The limits of the previous finite integrations will depend, as it was said
in the num$%
{{}^\circ}%
$.31., by the surfaces which are the boundaries of the body in the
configuration preceding the real one. The expression (7) will then be
adapted to represent the first, the second, the $n-th$ \ of the horizontal
series which are composing the quantity (4), by changing in it the
coordinates $a,b,c$ \ of the generic molecule $m_{i}$, that is giving to
these variables those suitable values which are needed for it to become one
after the others the molecules $m_{_{1}},m_{_{2}}....m_{n}$ ; and as the
said horizontal series are exactly $n$ (and also the terms of each of these
series are $n)$ the sum of all sums will be given to us by the finite
sextuple integral

\begin{equation}
\Sigma \delta a\Sigma \delta b\Sigma \delta c\Sigma \Delta f\Sigma \Delta
g\Sigma \Delta k\cdot \frac{1}{2}m_{i}m_{j}K\delta \rho .  \tag{9}
\label{CAPOVI-9ENG}
\end{equation}%
Let us recall now what we said in the last lines of the num$%
{{}^\circ}%
$.21., about the need of assigning the value $\sigma ^{3}$ to the letter $m$
which expresses the mass of the generic molecule: and as in the previous
integral there appears the product of two similar $m$ it will appear
manifest that this product must be replaced by the expression $\sigma
^{3}\cdot \sigma ^{3}$. Once we will have also recalled the analytical
theorem written in the equation (17) in the num$%
{{}^\circ}%
$.26.,

\subsubsection{Page 150}

\noindent theorem of which we will repeat here six times the application, we
will be ready to admit that the preceding finite sextuple integral is
transformed into the continuous sextuple integral

\begin{equation}
\int da\int db\int dc\int df\int dg\int dk\cdot\frac{1}{2}K\delta\rho 
\tag{10}  \label{CAPOVI-10ENG}
\end{equation}
with the addition of other terms, which then must be neglected. In it the
integration limits for the variables $f,g,k$ will depend on the surfaces
which bound the body in the antecedent configuration, and also on the
position of the molecule $m$, which is kept constant, that is by the
variables $a,b,c$ which after the first three will also vary in the same
domain.

72. Let us now spend some time developing some considerations about the
internal force $K$ which is exerted between one molecule and another
molecule, being either attractive or repulsive forces, which would have
acted both independently of the applied external forces and because of the
presence of these external forces. To assume that, as it was at first
suggested, it is a function $K(\rho )$ of the molecular distance only, it is
admissible only in the case of fluid bodies, as in that case the part of the
action due to the shape of the molecules is not present. In general (the
reader is invited to read once more what said in the num$%
{{}^\circ}%
$.54.) when the action due to the shape of molecules cannot be neglected,
there [in the expression for $K$] must appear also those cosines $%
a_{1},a_{2},a_{3},a_{4},ec.$which are fixing the directions of the edges or
axes of the two [involved] molecules relatively to the three orthogonal
axes, cosines whose values are changing from one molecule to the other and
therefore must result to be functions of the corresponding coordinates. It
could be very difficult to find how such functions have to be assigned (and
it is sufficient to this aim only to imagine that said directions could be
normal to curved surfaces having various and unknown nature for different
bodies): and beyond the ignorance about the internal structure of these
functions, it is not known how the $K$ depends on them. As a consequence the 
$K,$ if one wants to keep its most general possible use, must be a function
of the six coordinates, whose values are expressed by the equations (5),
(6): and we cannot presume to express its form, as we can only argue that it
has to be symmetric relatively to said six

\subsubsection{Page 151}

\noindent variables when taken three by three: i.e. that when changing the $%
x_{i},y_{i},z_{i}$ into the $x_{j},y_{j},z_{j},$ and these last into the
first ones the [function $K$ ] will remain the same. This is true because it
is known (as there is no reason for the contrary) that one half of $K$
expresses the action of the molecule $m_{j}$ on the molecule $m_{i}$, and
the other half of $K$ expresses the reciprocal action: it is possible to
assume that the roles of the two molecules are exchanged, and
notwithstanding this the analytical values must remain the same: this
observation leads us to conclude the stated property of the analytical
expression, as we mentioned also in the previous num$%
{{}^\circ}%
$. The impossibility of assigning the function $%
K(x_{i},y_{i},z_{i},x_{j},y_{j},z_{j})$ can be deduced also by means of
other arguments which I wish to omit: only I will remark how also from this
point of view the superiority of the methods which we have in our hands is
emerging: with them one can continue safely to proceed in the argumentation
notwithstanding an ignorance which cannot be defeated. We will add another
observation about the smallness of this molecular force $K$, by recalling
what we said about this subject at the end of the num$%
{{}^\circ}%
$.22. \ Similarly to what expounded in the num$%
{{}^\circ}%
$ 18 and following ones, it has to be regarded as an elementary force which
is so small that once considering the resultant of such forces on a single
physical point as acted by all the other molecules of the mas, we have still
a very small force of the same order of those forces $\sigma ^{3}X,\sigma
^{3}Y,\sigma ^{3}Z$ considered in the num$%
{{}^\circ}%
$.18. To this concept corresponds perfectly the scaling given by the factor $%
\sigma ^{6}$, which we will see to result from the sextuple integral (9) due
to the product $m_{i}m_{j}$ of the two elementary masses.

As a consequence of what was were said up to now we can, by adding up the
two integrals (1), (10), and by replacing the obtained sum in the first two
parts of the general equation (1) num$%
{{}^\circ}%
$.16., formulate the equation which includes the whole molecular mechanics.
Before doing so we will remark that it is convenient to introduce the
following definition

\begin{equation}
\Lambda =\frac{1}{4}\frac{K}{\rho }  \tag{11}  \label{CAPOVI-11ENG}
\end{equation}%
by means of which it will be possible to introduce the quantity $\Lambda
\delta \rho ^{2}$ instead of the quantity $\frac{1}{2}K\delta \rho $ in the
sextuple integral (10); and that

\subsubsection{Page 152}

\noindent inside this sextuple integral it is suitable to isolate the part
relative to the triple integral relative to the variables $f,g,k,$ placing
it under the same sign of triple integral with respect to the variables $%
a,b,c$ \ which includes the first part of the equation: which is manifestly
allowed. In this way the aforementioned general equation becomes

\begin{equation*}
\int da\int db\int dc\cdot\left\{ \left( X-\frac{d^{2}x}{dt^{2}}\right)
\delta x+\left( Y-\frac{d^{2}y}{dt^{2}}\right) \delta y+\left( Z-\frac{d^{2}z%
}{dt^{2}}\right) \delta z\right.
\end{equation*}

\begin{equation}
\left. +\int df\int dg\int dk\cdot \Lambda \delta \rho ^{2}\right\} +\Omega
=0  \tag{12}  \label{CAPOVI-12ENG}
\end{equation}%
where it is intended that (as mentioned at the beginning of the num$%
{{}^\circ}%
$.71.) it is included in the $\Omega $ the whole part which may be
introduced because of the forces applied to surfaces, lines or
well-determined points and also because of particular conditions which may
oblige some points to belong to some given curve or surface. This equation
(12) replaces the equation (1) of the num$%
{{}^\circ}%
$.44., or the equation (12) of the num$%
{{}^\circ}%
$.46$.$, and it is seen how it is expressed differently (while the remaining
parts are left the same) the part introduced by the reciprocal actions of
the molecules which in those equations were taken into account by means of
equations of conditions to hold in the whole body.

73. What remains to be done in order to deduce useful consequences from the
equation (12) is simply a calculation process. Once recalled the equation
(8), it is seen, transforming into series the functions in the brackets, so
that one has

\begin{align*}
\rho^{2} & =\left( f\frac{\text{ }dx}{da}+g\frac{dx}{db}+k\frac{dx}{dc}+%
\frac{f\text{ }^{2}}{2}\frac{d^{2}x}{da^{2}}+ec.\right) ^{2} \\
& +\left( f\frac{\text{ }dy}{da}+g\frac{dy}{db}+k\frac{dy}{dc}+\frac{f\text{ 
}^{2}}{2}\frac{d^{2}y}{da^{2}}+ec.\right) ^{2} \\
& +\left( f\frac{\text{ }dz}{da}+g\frac{dz}{db}+k\frac{dz}{dc}+\frac{f\text{ 
}^{2}}{2}\frac{d^{2}z}{da^{2}}+ec.\right) ^{2};
\end{align*}
and by calculating the squares and gathering the terms which have equal
coefficients:

\begin{align}
\rho^{2} & =f\text{ }^{2}t_{1}+g\text{ }^{2}t_{2}+k\text{ }%
^{2}t_{3}+2fgt_{4}+2fkt_{5}+2gkt_{6}  \notag \\
& +f\text{ }^{3}T_{1}+2f\text{ }^{2}gT_{2}+2f\text{ }^{2}kT_{3}+f\text{ }%
g^{2}T_{4}+ec.  \tag{13}  \label{CAPOVI-13ENG}
\end{align}

\subsubsection{Page\ 153}

\noindent in which expression the quantities $%
t_{1},t_{2},t_{3},t_{4},t_{5},t_{6}$ represent the six trimonials which are
alreay familiar to us, as we have adopted such denominations since the
equations (6) in the num$%
{{}^\circ}%
$.34.; and the quantities $T_{1},T_{2},T_{3},T_{4},ec.$ where the index goes
to infinity, represent trinomials of the same nature in which derivatives of
higher and higher order appear. In all these trinomials the last two terms
are always similar to the first but they differ in having the letters $y,z$
at the place of the letter $x.$ Those in which the second derivatives appear
are of two kinds. There are those which are composed with first order and
second order derivatives, and these are exactly 18 in number, which are
listed in the following formula: 
\begin{align}
& \frac{dx}{da}\frac{d^{2}x}{da^{2}}+\frac{dy}{da}\frac{d^{2}y}{da^{2}}+%
\frac{dz}{da}\frac{d^{2}z}{da^{2}}  \notag \\
& \frac{dx}{da}\frac{d^{2}x}{dadb}+\frac{dy}{da}\frac{d^{2}y}{dadb}+\frac{dz%
}{da}\frac{d^{2}z}{dadb}  \notag \\
& \frac{dx}{da}\frac{d^{2}x}{dadc}+\frac{dy}{da}\frac{d^{2}y}{dadc}+\frac{dz%
}{da}\frac{d^{2}z}{dadc}  \notag \\
& \frac{dx}{db}\frac{d^{2}x}{dadb}+\frac{dy}{db}\frac{d^{2}y}{dadb}+\frac{dz%
}{db}\frac{d^{2}z}{dadb}  \notag \\
& \frac{dx}{db}\frac{d^{2}x}{db^{2}}+\frac{dy}{db}\frac{d^{2}y}{db^{2}}+%
\frac{dz}{db}\frac{d^{2}z}{db^{2}}  \notag \\
& \frac{dx}{db}\frac{d^{2}x}{dbdc}+\frac{dy}{db}\frac{d^{2}y}{dbdc}+\frac{dz%
}{db}\frac{d^{2}z}{dbdc}  \notag \\
& \frac{dx}{dc}\frac{d^{2}x}{dadc}+\frac{dy}{dc}\frac{d^{2}y}{dadc}+\frac{dz%
}{dc}\frac{d^{2}z}{dadc}  \notag \\
& \frac{dx}{dc}\frac{d^{2}x}{dbdc}+\frac{dy}{dc}\frac{d^{2}y}{dbdc}+\frac{dz%
}{dc}\frac{d^{2}z}{dbdc}  \notag \\
& \frac{dx}{dc}\frac{d^{2}x}{dc^{2}}+\frac{dy}{dc}\frac{d^{2}y}{dc^{2}}+%
\frac{dz}{dc}\frac{d^{2}z}{dc^{2}}  \tag{14}  \label{CAPOVI-14ENG} \\
& \frac{dx}{db}\frac{d^{2}x}{da^{2}}+\frac{dy}{db}\frac{d^{2}y}{da^{2}}+%
\frac{dz}{db}\frac{d^{2}z}{da^{2}}  \notag \\
& \frac{dx}{da}\frac{d^{2}x}{db^{2}}+\frac{dy}{da}\frac{d^{2}y}{db^{2}}+%
\frac{dz}{da}\frac{d^{2}z}{db^{2}}  \notag \\
& \frac{dx}{da}\frac{d^{2}x}{dbdc}+\frac{dy}{da}\frac{d^{2}y}{dbdc}+\frac{dz%
}{da}\frac{d^{2}z}{dbdc}  \notag
\end{align}

\subsubsection{Page 154}

\begin{align*}
& \frac{dx}{dc}\frac{d^{2}x}{da^{2}}+\frac{dy}{dc}\frac{d^{2}y}{da^{2}}+%
\frac{dz}{dc}\frac{d^{2}z}{da^{2}} \\
& \frac{dx}{db}\frac{d^{2}x}{dadc}+\frac{dy}{db}\frac{d^{2}y}{dadc}+\frac {dz%
}{db}\frac{d^{2}z}{dadc} \\
& \frac{dx}{da}\frac{d^{2}x}{dc^{2}}+\frac{dy}{da}\frac{d^{2}y}{dc^{2}}+%
\frac{dz}{da}\frac{d^{2}z}{dc^{2}} \\
& \frac{dx}{dc}\frac{d^{2}x}{dadb}+\frac{dy}{dc}\frac{d^{2}y}{dadb}+\frac {dz%
}{dc}\frac{d^{2}z}{dadb} \\
& \frac{dx}{dc}\frac{d^{2}x}{db^{2}}+\frac{dy}{dc}\frac{d^{2}y}{db^{2}}+%
\frac{dz}{dc}\frac{d^{2}z}{db^{2}} \\
& \frac{dx}{db}\frac{d^{2}x}{dc^{2}}+\frac{dy}{db}\frac{d^{2}y}{dc^{2}}+%
\frac{dz}{db}\frac{d^{2}z}{dc^{2}}
\end{align*}

Then we have the trinomial constituted with second order derivatives only,
and these last are 21 in number, and are listed in the following formula

\begin{align}
& \left( \frac{d^{2}x}{da^{2}}\right) ^{2}+\left( \frac{d^{2}y}{da^{2}}%
\right) ^{2}+\left( \frac{d^{2}z}{da^{2}}\right) ^{2}  \notag \\
& \left( \frac{d^{2}x}{db^{2}}\right) ^{2}+\left( \frac{d^{2}y}{db^{2}}%
\right) ^{2}+\left( \frac{d^{2}z}{db^{2}}\right) ^{2}  \notag \\
& \left( \frac{d^{2}x}{dc^{2}}\right) ^{2}+\left( \frac{d^{2}y}{dc^{2}}%
\right) ^{2}+\left( \frac{d^{2}z}{dc^{2}}\right) ^{2}  \notag \\
& \left( \frac{d^{2}x}{dadb}\right) ^{2}+\left( \frac{d^{2}y}{dadb}\right)
^{2}+\left( \frac{d^{2}z}{dadb}\right) ^{2}  \notag \\
& \left( \frac{d^{2}x}{dadc}\right) ^{2}+\left( \frac{d^{2}y}{dadc}\right)
^{2}+\left( \frac{d^{2}z}{dadc}\right) ^{2}  \notag \\
& \left( \frac{d^{2}x}{dbdc}\right) ^{2}+\left( \frac{d^{2}y}{dbdc}\right)
^{2}+\left( \frac{d^{2}z}{dbdc}\right) ^{2}  \notag \\
& \frac{d^{2}x}{da^{2}}\frac{d^{2}x}{db^{2}}+\frac{d^{2}y}{da^{2}}\frac {%
d^{2}y}{db^{2}}+\frac{d^{2}z}{da^{2}}\frac{d^{2}z}{db^{2}}  \notag \\
& \frac{d^{2}x}{da^{2}}\frac{d^{2}x}{dc^{2}}+\frac{d^{2}y}{da^{2}}\frac {%
d^{2}y}{dc^{2}}+\frac{d^{2}z}{da^{2}}\frac{d^{2}z}{dc^{2}}  \notag \\
& \frac{d^{2}x}{db^{2}}\frac{d^{2}x}{dc^{2}}+\frac{d^{2}y}{db^{2}}\frac {%
d^{2}y}{dc^{2}}+\frac{d^{2}z}{db^{2}}\frac{d^{2}z}{dc^{2}}  \notag \\
& \frac{d^{2}x}{da^{2}}\frac{d^{2}x}{dadb}+\frac{d^{2}y}{da^{2}}\frac{d^{2}y%
}{dadb}+\frac{d^{2}z}{da^{2}}\frac{d^{2}z}{dadb}  \notag \\
& \frac{d^{2}x}{da^{2}}\frac{d^{2}x}{dadc}+\frac{d^{2}y}{da^{2}}\frac{d^{2}y%
}{dadc}+\frac{d^{2}z}{da^{2}}\frac{d^{2}z}{dadc}  \tag{15}
\label{CAPOVI-15ENG}
\end{align}

\subsubsection{Page\ 155}

\begin{align*}
& \frac{d^{2}x}{da^{2}}\frac{d^{2}x}{dbdc}+\frac{d^{2}y}{da^{2}}\frac{d^{2}y%
}{dbdc}+\frac{d^{2}z}{da^{2}}\frac{d^{2}z}{dbdc} \\
& \frac{d^{2}x}{db^{2}}\frac{d^{2}x}{dadb}+\frac{d^{2}y}{db^{2}}\frac{d^{2}y%
}{dadb}+\frac{d^{2}z}{db^{2}}\frac{d^{2}z}{dadb} \\
& \frac{d^{2}x}{db^{2}}\frac{d^{2}x}{dadc}+\frac{d^{2}y}{db^{2}}\frac{d^{2}y%
}{dadc}+\frac{d^{2}z}{db^{2}}\frac{d^{2}z}{dadc} \\
& \frac{d^{2}x}{db^{2}}\frac{d^{2}x}{dbdc}+\frac{d^{2}y}{db^{2}}\frac{d^{2}y%
}{dbdc}+\frac{d^{2}z}{db^{2}}\frac{d^{2}z}{dbdc} \\
& \frac{d^{2}x}{dc^{2}}\frac{d^{2}x}{dadb}+\frac{d^{2}y}{dc^{2}}\frac{d^{2}y%
}{dadb}+\frac{d^{2}z}{dc^{2}}\frac{d^{2}z}{dadb} \\
& \frac{d^{2}x}{dc^{2}}\frac{d^{2}x}{dadc}+\frac{d^{2}y}{dc^{2}}\frac{d^{2}y%
}{dadc}+\frac{d^{2}z}{dc^{2}}\frac{d^{2}z}{dadc} \\
& \frac{d^{2}x}{dc^{2}}\frac{d^{2}x}{dbdc}+\frac{d^{2}y}{dc^{2}}\frac{d^{2}y%
}{dbdc}+\frac{d^{2}z}{dc^{2}}\frac{d^{2}z}{dbdc} \\
& \frac{d^{2}x}{dadb}\frac{d^{2}x}{dadc}+\frac{d^{2}y}{dadb}\frac{d^{2}y}{%
dadc}+\frac{d^{2}z}{dadb}\frac{d^{2}z}{dadc} \\
& \frac{d^{2}x}{dadb}\frac{d^{2}x}{dbdc}+\frac{d^{2}y}{dadb}\frac{d^{2}y}{%
dbdc}+\frac{d^{2}z}{dadb}\frac{d^{2}z}{dbdc} \\
& \frac{d^{2}x}{dadc}\frac{d^{2}x}{dbdc}+\frac{d^{2}y}{dadc}\frac{d^{2}y}{%
dbdc}+\frac{d^{2}z}{dadc}\frac{d^{2}z}{dbdc}.
\end{align*}
The trinomials with third order derivatives are of three kinds: there are
those constituted by derivatives of first and third order, and one can count
30 of them: there are those constituted by derivatives of second and third
order, and they are 60 in number: and there are those which contain only
third order derivatives and they are 55 in number. I am not writing them, as
everybody who is given the needed patience can easily calculate them by
himself, as it can be also done for those trinomials containing derivatives
of higher order.

When using the equation (13) to deduce the value of the variation $\delta
\rho ^{2}$ , it is clear that the characteristic $\delta $ will need to be
applied only to the trinomials we have discussed up to now, so that we will
have:

\begin{align}
\delta\rho^{2} & =f\text{ }^{2}\delta t_{1}+g\text{ }^{2}\delta t_{2}+k\text{
}^{2}\delta t_{3}+2fg\text{ }\delta t_{4}+2fk\text{ }\delta t_{5}+2gk\text{ }%
\delta t_{6}  \notag \\
& +f\text{ }^{3}\delta T_{1}+2f\text{ }^{2}g\delta T_{2}+2f\text{ }%
^{2}k\delta T_{3}+f\text{ }g^{2}\delta T_{4}+ec.  \tag{16}
\label{CAPOVI-16ENG}
\end{align}

\subsubsection{Page 156}

\noindent Indeed the coefficients $f$ $^{2},g$ $^{2},k$ $^{2},2fg$, etc. are
always of the same form as the functions giving the variables $x,y,z$ in
terms of the variables $a,b,c,$ and therefore cannot be affected by that
operation whose aim is simply to change the form of these functions. Vice
versa, by multiplying the previous equation (16) times $\Lambda $ and then
integrating with respect to the variables $f,g,k$ \ in order to deduce from
such calculation the value to be given to the forth term under the triple
integral of the equation (12), such an operation is affecting only the
quantities $\Lambda f$ $^{2},\Lambda g$ $^{2}$, etc. and the variations $%
\delta t_{1},\delta t_{2},\delta t_{3}....\delta T_{1},\delta T_{2},ec.$
cannot be affected by it, as the trinomials $%
t_{1},t_{2},t_{3}....T_{1},T_{2},ec.$ (one has to consider carefully which
is their origin) do not contain the variables $f,g,k$ : therefore such
variations result to be constant factors, times which are to be multiplied
the integrals to be calculated in the subsequent terms of the series. After
these considerations it is manifest the truth of the equation:

\begin{equation}
\int df\int dg\int dk\cdot\Lambda\delta\rho^{2}=  \tag{17}
\label{CAPOVI-17ENG}
\end{equation}

\begin{align*}
& \left( 1\right) \text{ }\delta t_{1}+\left( 2\right) \text{ }\delta
t_{2}+\left( 3\right) \text{ }\delta t_{3}+\left( 4\right) \text{ }\delta
t_{4}+\left( 5\right) \text{ }\delta t_{5}+\left( 6\right) \text{ }\delta
t_{6} \\
& +\left( 7\right) \text{ }\delta T_{1}+\left( 8\right) \text{ }\delta
T_{2}+\left( 9\right) \text{ }\delta T_{3}+\left( 10\right) \text{ }\delta
T_{4}+ec.
\end{align*}%
where the coefficients (1), (2), etc. indicated by means of numbers in
between brackets, must be regarded to be each a function of the variables $%
a,b,c$ as obtained after having performed the said definite integrals. This
is the equivalent quantity which should be introduced in the equation (12)
at the place of the forth term under the triplicate integral.

74. A new proposition, to which the reader should pay much attention, is
that all the trinomials $T_{1},T_{2},T_{3},etc.$ where the index goes to
infinity , which appear in the previous equation (17), can be expresses by
means of the only first six $t_{1},t_{2},t_{3},$ $t_{4},t_{5},t_{6,}$ and of
their derivatives with respect to the variables $a,b,c$ of all orders. I
started to suspect this analytical truth because of the necessary
correspondence which must hold between the results which are obtained with
the way considered in this Capo and those results obtained with the way
considered in the Capos III and IV. I have then verified the stated property
for 39 terms of the previous

\subsubsection{Page\ 157}

\noindent series (17), beyond the first six, that is for all trinomials
written in the tables (14), and (15), and after these calculations I
abandoned myself to the analogy: and this will be sooner or later
unavoidable, because our series is infinite and it will be impossible to
check all its terms. Now I will say how I performed the stated verification
and the importance of the conclusions will justify the lengthinesses of the
calculations, which, except for the prolixity, do not present any
difficulty. Checking the values of the variables $t_{1},t_{2},t_{3},$ $%
t_{4},t_{5},t_{6}$ ( equations (6) num$.%
{{}^\circ}%
$34.) it is immediate to recognize that the first nine trinomials of the
table (14) respectively have the values:

\begin{align*}
& \frac{1}{2}\frac{dt_{1}}{da};\text{ \ \ \ }\frac{1}{2}\frac{dt_{1}}{db};%
\text{ \ \ \ }\frac{1}{2}\frac{dt_{1}}{dc} \\
& \frac{1}{2}\frac{dt_{2}}{da};\text{ \ \ \ }\frac{1}{2}\frac{dt_{2}}{db};%
\text{ \ \ \ }\frac{1}{2}\frac{dt_{2}}{dc} \\
& \frac{1}{2}\frac{dt_{3}}{da};\text{ \ \ \ }\frac{1}{2}\frac{dt_{3}}{db};%
\text{ \ \ \ }\frac{1}{2}\frac{dt_{3}}{dc}.
\end{align*}

It is then possible to find that (and this can be verified by simple
substitution of known values) the trinomials labelled with the number ten,
eleven, thirteen, fifteen, seventeen and eighteen are equivalent
respectively to the following binomials:

\begin{align*}
& \frac{dt_{4}}{da}-\frac{1}{2}\frac{dt_{1}}{db};\text{ \ \ \ }\frac{dt_{4}}{%
db}-\frac{1}{2}\frac{dt_{2}}{da};\text{ \ \ \ }\frac{dt_{5}}{da}-\frac{1}{2}%
\frac{dt_{1}}{dc} \\
& \frac{dt_{5}}{dc}-\frac{1}{2}\frac{dt_{3}}{da};\text{ \ \ \ }\frac{dt_{6}}{%
db}-\frac{1}{2}\frac{dt_{2}}{dc};\text{ \ \ \ }\frac{dt_{6}}{dc}-\frac{1}{2}%
\frac{dt_{3}}{db}.
\end{align*}
And that the trinomials labelled with the numbers twelve, fourteen and
sixteen have respectively these other values:

\begin{equation*}
\frac{1}{2}\frac{dt_{4}}{dc}+\frac{1}{2}\frac{dt_{5}}{db}-\frac{1}{2}\frac{%
dt_{6}}{da};\text{ \ \ \ }\frac{1}{2}\frac{dt_{4}}{dc}-\frac{1}{2}\frac{%
dt_{5}}{db}+\frac{1}{2}\frac{dt_{6}}{da};\text{ \ \ }-\frac{1}{2}\frac{dt_{4}%
}{dc}+\frac{1}{2}\frac{dt_{5}}{db}+\frac{1}{2}\frac{dt_{6}}{da}.
\end{equation*}
In this way the stated proposition is proven relatively to the first 18
trinomials.

Let us now imagine to have formed 18 equations having each as left hand
sides the trinomials of the table (14) and as right hand side respectively
the values to which we have proven they are equal. Among these equations let
us consider the

\subsubsection{Page 158}

\noindent first, the tenth and the thirteenth ones and let us multiply them
orderly times $l_{1},m_{1},n_{1}$, and then let us sum the obtained results:
then let us multiply the same equations again times $l_{2},m_{2},n_{2}$, and
again sum the obtained results : finally let us multiply the same three
equations times $l_{3},m_{3},n_{3}$, and again sum the obtained results : by
considering the nine equations of the num.$%
{{}^\circ}%
$14. labelled (28) , we will manage to get an expression for the values of
the three second order derivatives $\frac{d^{2}x}{da^{2}},\frac{d^{2}y}{%
da^{2}},\frac{d^{2}z}{da^{2}}$. Following the same procedure, suitably
choosing among the aforementioned 18 equations, we will determine the values
of the other second order derivatives and we will get

\begin{align*}
2H\text{ }\frac{d^{2}x}{da^{2}} & =l_{1}\frac{dt_{1}}{da}+m_{1}\left( 2\frac{%
dt_{4}}{da}-\frac{dt_{1}}{db}\right) +n_{1}\left( 2\frac{dt_{5}}{da}-\frac{%
dt_{1}}{dc}\right) \\
2H\text{ }\frac{d^{2}y}{da^{2}} & =l_{2}\frac{dt_{1}}{da}+m_{2}\left( 2\frac{%
dt_{4}}{da}-\frac{dt_{1}}{db}\right) +n_{2}\left( 2\frac{dt_{5}}{da}-\frac{%
dt_{1}}{dc}\right) \\
2H\text{ }\frac{d^{2}z}{da^{2}} & =l_{3}\frac{dt_{1}}{da}+m_{3}\left( 2\frac{%
dt_{4}}{da}-\frac{dt_{1}}{db}\right) +n_{3}\left( 2\frac{dt_{5}}{da}-\frac{%
dt_{1}}{dc}\right) \\
2H\text{ }\frac{d^{2}x}{db^{2}} & =l_{1}\left( 2\frac{dt_{4}}{db}-\frac{%
dt_{2}}{da}\right) +m_{1}\frac{dt_{2}}{db}+n_{1}\left( 2\frac{dt_{6}}{db}-%
\frac{dt_{2}}{dc}\right) \\
2H\text{ }\frac{d^{2}y}{db^{2}} & =l_{2}\left( 2\frac{dt_{4}}{db}-\frac{%
dt_{2}}{da}\right) +m_{3}\frac{dt_{2}}{db}+n_{2}\left( 2\frac{dt_{6}}{db}-%
\frac{dt_{2}}{dc}\right) \\
2H\text{ }\frac{d^{2}z}{db^{2}} & =l_{3}\left( 2\frac{dt_{4}}{db}-\frac{%
dt_{2}}{da}\right) +m_{3}\frac{dt_{2}}{db}+n_{3}\left( 2\frac{dt_{6}}{db}-%
\frac{dt_{2}}{dc}\right) \\
2H\text{ }\frac{d^{2}x}{dc^{2}} & =l_{1}\left( 2\frac{dt_{5}}{dc}-\frac{%
dt_{3}}{da}\right) +m_{1}\left( 2\frac{dt_{6}}{dc}-\frac{dt_{3}}{db}\right)
+n_{1}\frac{dt_{3}}{dc} \\
2H\text{ }\frac{d^{2}y}{dc^{2}} & =l_{2}\left( 2\frac{dt_{5}}{dc}-\frac{%
dt_{3}}{da}\right) +m_{2}\left( 2\frac{dt_{6}}{dc}-\frac{dt_{3}}{db}\right)
+n_{2}\frac{dt_{3}}{dc} \\
2H\text{ }\frac{d^{2}z}{dc^{2}} & =l_{3}\left( 2\frac{dt_{5}}{dc}-\frac{%
dt_{3}}{da}\right) +m_{3}\left( 2\frac{dt_{6}}{dc}-\frac{dt_{3}}{db}\right)
+n_{3}\frac{dt_{3}}{dc} \\
2H\text{ }\frac{d^{2}x}{dadb} & =l_{1}\frac{dt_{1}}{db}+m_{1}\frac{dt_{2}}{da%
}+n_{1}\left( \frac{dt_{5}}{db}+\frac{dt_{6}}{da}-\frac{dt_{4}}{dc}\right) \\
2H\text{ }\frac{d^{2}y}{dadb} & =l_{2}\frac{dt_{1}}{db}+m_{2}\frac{dt_{2}}{da%
}+n_{2}\left( \frac{dt_{5}}{db}+\frac{dt_{6}}{da}-\frac{dt_{4}}{dc}\right) \\
2H\text{ }\frac{d^{2}z}{dadb} & =l_{3}\frac{dt_{1}}{db}+m_{3}\frac{dt_{2}}{da%
}+n_{3}\left( \frac{dt_{5}}{db}+\frac{dt_{6}}{da}-\frac{dt_{4}}{dc}\right) \\
2H\text{ }\frac{d^{2}x}{dadc} & =l_{1}\frac{dt_{1}}{dc}+m_{1}\left( \frac{%
dt_{4}}{dc}+\frac{dt_{6}}{da}-\frac{dt_{5}}{db}\right) +n_{1}\frac{dt_{3}}{da%
}
\end{align*}

\subsubsection{Page 159}

\begin{align*}
2H\text{ }\frac{d^{2}y}{dadc}& =l_{2}\frac{dt_{1}}{dc}+m_{2}\left( \frac{%
dt_{4}}{dc}+\frac{dt_{6}}{da}-\frac{dt_{6}}{db}\right) +n_{2}\frac{dt_{3}}{da%
} \\
2H\text{ }\frac{d^{2}z}{dadc}& =l_{3}\frac{dt_{1}}{dc}+m_{3}\left( \frac{%
dt_{4}}{dc}+\frac{dt_{6}}{da}-\frac{dt_{5}}{db}\right) +n_{3}\frac{dt_{3}}{da%
} \\
2H\text{ }\frac{d^{2}x}{dbdc}& =l_{1}\left( \frac{dt_{4}}{dc}+\frac{dt_{5}}{%
db}-\frac{dt_{6}}{da}\right) +m_{1}\frac{dt_{2}}{dc}+n_{1}\frac{dt_{3}}{db}
\\
2H\text{ }\frac{d^{2}y}{dbdc}& =l_{2}\left( \frac{dt_{4}}{dc}+\frac{dt_{5}}{%
db}-\frac{dt_{6}}{da}\right) +m_{2}\frac{dt_{2}}{dc}+n_{2}\frac{dt_{3}}{db}
\\
2H\text{ }\frac{d^{2}z}{dbdc}& =l_{3}\left( \frac{dt_{4}}{dc}+\frac{dt_{5}}{%
db}-\frac{dt_{6}}{da}\right) +m_{3}\frac{dt_{2}}{dc}+n_{3}\frac{dt_{3}}{db}.
\end{align*}%
Now, by means of these values, let us look for the values of the trinomials
of the table (15). Recalling the equations (31), (33), (34) in the num.$%
{{}^\circ}%
$ 67. we will see that these values result to be constituted uniquely by the
variables $t_{1},t_{2}....t_{6}$ and by their first order derivatives, and
this is exactly what we wanted to prove. For instance the value of the first
trinomial

\begin{equation*}
\left( \frac{d^{2}x}{da^{2}}\right) ^{2}+\left( \frac{d^{2}y}{da^{2}}\right)
^{2}+\left( \frac{d^{2}z}{da^{2}}\right) ^{2}
\end{equation*}
can be proven to be equal to a fraction whose numerator is

\begin{align*}
& \left( t_{_{2}}t_{_{3}}-t_{_{6}}^{^{2}}\right) \left( \frac{dt_{_{1}}}{da}%
\right) ^{2}+\left( t_{_{1}}t_{_{3}}-t_{_{5}}^{^{2}}\right) \left( 2\frac{%
dt_{4}}{da}-\frac{dt_{1}}{db}\right) ^{2}+\left(
t_{_{1}}t_{_{2}}-t_{_{4}}^{^{2}}\right) \left( 2\frac{dt_{5}}{da}-\frac{%
dt_{1}}{dc}\right) ^{2} \\
& +2\left( t_{_{5}}t_{_{6}}-t_{_{3}}t_{_{4}}\right) \frac{dt_{_{1}}}{da}%
\left( 2\frac{dt_{4}}{da}-\frac{dt_{1}}{db}\right) +2\left(
t_{_{4}}t_{_{6}}-t_{_{2}}t_{_{5}}\right) \frac{dt_{_{1}}}{da}\left( 2\frac{%
dt_{6}}{da}-\frac{dt_{1}}{dc}\right) \\
& +2\left( t_{_{4}}t_{_{5}}-t_{_{1}}t_{_{6}}\right) \left( 2\frac{dt_{4}}{da}%
-\frac{dt_{1}}{db}\right) \left( 2\frac{dt_{5}}{da}-\frac{dt_{1}}{dc}\right)
\end{align*}
and whose denominator is:

\begin{equation*}
4\left(
t_{_{1}}t_{_{2}}t_{_{3}}+2t_{_{4}}t_{_{5}}t_{_{6}}-t_{_{1}}t_{_{6}}^{^{2}}-t_{_{2}}t_{_{5}}^{^{2}}-t_{_{3}}t_{_{4}}^{^{2}}\right) .
\end{equation*}%
Similarly forms can be found for the values of the other twenty trinomials
of the table (15) : therefore it was not exaggerated to affirm that the
stated analytical property has been actually verified for 39 trinomials.

75. Once the proposition of the previous num. has been admitted, it is
manifest that the equation (17) can assume the following other form

\subsubsection{Page 160}

\begin{equation}
\int df\int dg\int dk\cdot\Lambda\delta\rho^{2}=  \tag{18}
\label{CAPOVI-18ENG}
\end{equation}

\begin{align*}
& \left( \alpha \right) \text{ }\delta t_{1}+\left( \beta \right) \text{ }%
\delta t_{2}+\left( \gamma \right) \text{ }\delta t_{3}+....+\left( \epsilon
\right) \text{ }\frac{\delta dt_{_{1}}}{da}+\left( \zeta \right) \text{ }%
\frac{\delta dt_{_{1}}}{db}+\left( \eta \right) \text{ }\frac{\delta
dt_{_{1}}}{dc} \\
& +\left( \vartheta \right) \text{ }\frac{\delta dt_{_{2}}}{da}+....+\left(
\lambda \right) \text{ }\frac{\delta d^{2}t_{_{1}}}{da^{2}}+\left( \mu
\right) \text{ }\frac{\delta d^{2}t_{_{1}}}{dadb}+....+\left( \xi \right) 
\text{ }\frac{\delta d^{2}t_{_{2}}}{da^{2}}+\left( o\right) \frac{\delta
d^{2}t_{_{2}}}{dadb}+ec.
\end{align*}%
in which the coefficients $\left( \alpha \right) ,\left( \beta \right) $ $%
....\left( \epsilon \right) ....\left( \lambda \right) ....ec.$ \ are
suitable quantities given in terms of the coefficients $\left( 1\right)
,\left( 2\right) ....\left( 7\right) ,\left( 8\right) ....$ of the equation
(17): they depend on the quantities $t_{1},t_{2}....t_{6}$, and on all order
derivatives of these trinomials with respect to the variables $a,b,c$ . Then
the variations $\delta t_{1},$ $\delta t_{2}...$.(with the index varying up
to infinity) and the variations of all their derivatives of all orders $%
\frac{\delta dt_{_{1}}}{da},\frac{\delta dt_{_{1}}}{db},ec.$ appear in the
(18) only linearly. Now it is a fundamental principle of the calculus of
variations (and we used it also in this Memoir in the num.$%
{{}^\circ}%
$ 36. and elsewhere) that one series as the previous one, where the
variations of some quantities and the variations of their derivatives with
respect to the fundamental variables $a,b,c$ appear linearly can be always
be transformed into one expression which containes the first quantities
without any sign of derivation, with the addition of other terms which are
exact derivatives with respect to one of the three simple independent
variables. As a consequence of such a principle, the expression which
follows to the equation (18) can be given

\begin{equation}
\int df\int dg\int dk\cdot\Lambda\delta\rho^{2}=  \tag{19}
\label{CAPOVI-19ENG}
\end{equation}

\begin{align*}
& A\delta t_{1}+B\text{ }\delta t_{2}+C\text{ }\delta t_{3}+D\delta
t_{4}+E\delta t_{5}+F\text{ }\delta t_{6} \\
& \text{ \ \ \ \ \ \ \ \ \ \ \ \ \ \ \ }+\frac{d\Delta }{da}+\frac{d\Theta }{%
db}+\frac{d\Upsilon }{dc}.
\end{align*}%
The values of the six coefficients $A,B,C,D,E,F$ are series constructed with
the coefficients $\left( \alpha \right) ,\left( \beta \right) ,\left( \gamma
\right) ....\left( \epsilon \right) ,\left( \zeta \right) ....\left( \lambda
\right) ,ec.$ of the equation (18) which appear linearly, with alternating
signs and affected by derivations of higher and higher order when we move
ahead in the terms of said series: the quantities $\Delta ,\Theta ,\Upsilon $
\ are series of the same form of the terms which are transformed, in which
the coefficients of the variations have a composition

\subsubsection{Page\ 161}

\noindent similar to the one which we have described for the six
coefficients $A,B,C,D,E,F$ .

Once -instead of the quantity under the integral sign in the left hand side
of the equation (12)- one introduces the quantity which is on the right hand
side of the equation (19), it is clear to everybody that an integration is
possible for each of the last three addends appearing in it and that as a
consequence these terms only give quantities which supply boundary
conditions. What remains under the triple integral is the only sestinomial
which is absolutely similar to the sestinomial already used in the equation
(10) num.$%
{{}^\circ}%
$ 35. for rigid systems. Therefore after having remarked the aforementioned
similarity the analytical procedure to be used here will result perfectly
equal to the one used in the num.$%
{{}^\circ}%
$ 35, procedure which led to the equations (26), (29) in the num.$%
{{}^\circ}%
$ 38 and it will become possible the demonstration of the extension of the
said equations to every kind of bodies which do not respect the constraint
of rigidity, as it was mentioned at the end of the num.$%
{{}^\circ}%
$ 38. It will also be visible the coincidence of the obtained results with
those which are expressed in the equations (23) of the num.$%
{{}^\circ}%
$50. which hold for every kind of systems and which were shown in the Capo
IV by means of those intermediate coordinates $p,q,r$, whose consideration,
when using the approach used in this Capo, will not be needed.

76. The preceding analysis allows for many useful considerations. First of
all I will mention those which are needed to clarify the doubts to which we
already tried to answer in the num.$%
{{}^\circ}%
$ 63., when we promised to add more explanations exactly in this later
number. Starting from the antecedent ideal configuration of the molecules
described by means of the coordinates $a,b,c$, and arriving at the real
configuration, we intend this second configuration as described by means of
two different reference triples of orthogonal axes, those of the variables $%
p,q,r$, and of the variables $x,y,z$. \ To express the real configuration by
means of the variables $p,q,r$ we need simply to copy the preceding analysis
writing everywhere $p,q,r$ where previously we had written $x,y,z$. \ Now if
we want to pass from the coordinates $p,q,r$ to the coordinates $x,y,z$, we
will observe that in the case of the fluids, as it was said in the num.$%
{{}^\circ}%
$ 72., the internal force $K$, or equivalently the quantity $\Lambda$,

\subsubsection{Page 162}

\noindent is uniquely a function of the molecular distance $\rho $ ; and if
one considers carefully the [algebraic and differential] form of the
coefficients $(1),(2),(3)......$ in the equation (17), and that of the
coefficients $(\alpha ),(\beta ),(\gamma )......$ of the equation (18), and
finally that of the six coefficients $A,B,C,D,E,F$ \ in the equation (19),
we will be easily persuaded that in these last six quantities, when
introducing the first reference frame relative to the variables $p,q,r$,
will contain such variables $p,q,r$ only because these last appear in the
radical $\rho $ (equation (8)), and in the six trinomials $%
t_{1},t_{2},....t_{6}$ , having written everywhere the letters $p,q,r$ at
the place of the letters $x,y,z$. Therefore the six quantities $A,B,C,D,E,F$
will enjoy the analytical property which we have discussed many times, which
consists in changing when undergoing the substitution of the values (31) num.%
$%
{{}^\circ}%
$ 40. in quantities which depends in an equal way on the variables $x,y,z,$
as any trace of the nine angular quantities $\alpha _{1,}\beta _{1,}\gamma
_{1,}\alpha _{2,}ec.,$ disappears because such property is verified by the
radical $\rho $ , and by each of the six trinomials $t_{1},t_{2}....t_{6}$ .
Now in this case one has only to perform the indicated operations and he
will be persuaded that the situation is exactly the one which we have
described, once the equations among the angular quantities introduced in the
num.$%
{{}^\circ}%
$ 33 are recalled. Here is one of those results forecasting which, at the
beginning of the present Capo, we stated preliminarily that the two
conceived ways for calculating the internal constraints among the molecules
were illustrating each other.

Vice versa the procedure used in the Capo IV. makes possible the treatment
of the quantities at the boundaries which would be very difficult with the
present method. We have seen there (num.$%
{{}^\circ}%
$ 52.) as in the boundary conditions the same six quantities appear which
express the effect of internal forces in the three [bulk] equations valid in
all points of the body mass: here, however, the boundary conditions would be
complicated because of the role played by those other quantities $\Delta
,\Theta ,\Upsilon $ which appear in the last terms of the equation (19). It
is convenient to say that the whole part introduced by such terms in the
limit quantities is playing only an apparent role: and this analytical fact
is often encountered in the calculus of variations. In the questions related
to the calculus of variations which can be referred to formulas involving
definite triple integrals, if one has an equation of condition $L=0$ , and
considers the differential equations

\subsubsection{Page\ 163}

\begin{equation*}
\frac{dL}{da}=0,\frac{dL}{db}=0,\frac{dL}{dc}=0;\frac{d^{^{2}}L}{da^{^{2}}}%
=0, ec.
\end{equation*}%
of all orders, and if he treats them as if they were many new equations of
conditions, multiplying their variations times undetermined coefficients and
introducing the obtained products in the equations for maximum or minimum
under the integral sign and operating the usual transformations he will add
(to the quantity which would have been left without doing the further
operations which we have said) a quantity of the same nature of the
trinomial which is at the end of the equation (19). In such a case one
understands that the novelty of the appearance of quantities at the
boundaries can be only aparent, because the aforesaid differentiated
equations do not have any new meaning which was not already expresses by the
equation $L=0$ \ alone.

However the requirement that the said boundary quantity must vanish, as here
we are not dealing with a result concerning equations of conditions, which
cannot be used more than really necessary, but simply with a necessary
consequence of the comparison of the two presented methods, [that
requirement] could be a tool which can lead us to the discovery of new
truths, which could lead us, for instance, to more detailed investigations
about the nature of molecular interactions. I limit myself to a simple hint:
but concerning the molecular action I cannot avoid to remark something which
seems to me to be not irrelevant. The reader will have discovered that in
the previous analysis to come to the most important conclusions, that is to
the equations which hold for all points of the mass, it was not needed to
use the hypotheses accepted by modern Geometers, and by myself in the \S .V
in the Memoir which was published in the Tomo\ XXI of these Atti. It was
said that the molecular action must be appreciable at distances which cannot
be sensed and not appreciable for distances which can be sensed: this
statement may be true, but whatever it will be, one can ignore, in our
discussion, such a hypothesis: the integrals with respect to the variables $%
s,g,k$ in the equations (17), (18), (19) can be regarded to be extended up
at the boundaries of the body, and not only for very small distances, what
is done with some effort.

I will conclude the chapter by imagining an objection. Somebody could say
that: "you have here calculated the internal actions between

\subsubsection{Page 164}

\noindent one molecule and another molecule, and we can agree about this
point: but why you have not done the same also in the Capo IV, where you
left out all the second part of the general equation (1) num.$%
{{}^\circ}%
$ 16 ? Such omission can nullify all the deductions obtained with those
analytical procedures". I answer as follows:\ "I have not taken into account
then all the internal actions between pairs of molecules because at their
place I considered the six equations of condition which were treated there:
indeed considering also said actions would have meant to include their
effects two times. Also when one considers the motion of rigid bodies the
action between each pair of molecules is effective, but everybody who knows
the spirit of Analytical Mechanics will be persuaded that all of these
actions were accounted for by the six equations (8) in the num.$%
{{}^\circ}%
$ 34. which are valid for all physical points: on the contrary if one doubt
could be raised it could concern the fact that the six equations of
conditions were redundant, and for this reason the considerations in the
number 39 became necessary.

The same has to be said in the general case of the equations of conditions
(14) in the num.$%
{{}^\circ}%
$ 47., which also hold for all points of the mass: the only fact of
considering and calculating them is equivalent to account for all internal
forces, although it is not explicitly detailed how this occurs. To somebody
it could seem -in some aspects- more persuasive the method used in this
Capo, as it allows us to have an explanation about the way in which the
molecular actions are acting: however I believe that those who think in the
best way will prefer the method presented in the Capo IV, which is a method
more direct and more powerful, because it is based on that geometrical
principle which assumes the indifference of the choice of the reference
frames with respect to the system, principle which we will need also in the
next Capo, and which contains the reason of so many mechanical truths.
Moreover I see that it is possible -by following the flow of reasonings
expounded in the Capo IV- to explain the true meaning of that other
lagrangian principle, which at the beginning of the same Capo we said to be
too abstract, and to firmly establish its extension and the effective way of
using it. Aforementioned explication, however, would not be maybe of great
utility, and I state this because I believe that all advantages to which
Lagrange was aiming by postulating that principle could be obtained more
directly and naturally by using the procedure described in the same Capo IV.

\subsection{Appendix B Eulogy in memoriam of Vincenzo Brunacci by Gabrio
Piola.}

It is extremely painful for us to announce in this document the death of a
truly great man, who, as during his life was a glory for Italy, now moves,
because of his loss, everybody inconsolably to tears. One of the most
eminent mathematicians of our time, the illustrious professor Vincenzo
Brunacci of Pavia, suffering for many years because of a painful disease, on
the day 16th of the current month of June was attacked by those very strong
convulsions which were consequence of it, ceased to live surrounded by
friendship and religion. To recall the merits of the deceased person during
the time while still everybody cries on his tomb it is really a way for
increasing painful laments, even if it will be simply a meagre tribute of
praise which will be written by our pen. Our intention is not that of
presenting a formal elogium; this kind of encomium soon will be heard in the
most erudite academies and in all palaces of sciences.

Vincenzo Brunacci was born in the fatherland of Galileo the day 3rd of March
1768, his father first name being Ignazio Maria and his mother being named
Elisabetta Danielli.

It seemed as if the Spirit of Italy who was in great sufferance because in
that time the most brilliant star of all mathematical sciences, the
illustrious Lagrangia, had left the Nation, that Spirit wanted to have the
rise of another star, which being born on the banks of the river Arno, was
bound to become the successor of the first one.

This consideration is presenting itself even more spontaneous by when we
will remark that Brunacci was the first admirer in Italy of the luminous
Lagrangian doctrines, the scientist who diffused and supported them, the
scientist who in his studies was always a very creative innovator in their
applications.

His first Maestri were two famous Italians, Father Canovai and the great
geometer Pietro Paoli. Although in his first youth he was diverted by other
studies, which were opposed to his natural inclinations and from which,
because of due respect he could not subtract himself, he still was able to
cultivate at the same time those studies for which he had been born.

Very soon he was pupil only of the classic textbooks and of himself. Very
soon, as he did not allow to any man to see his genius while being born or
in his first childhood, in the "opuscolo analitico", printed in Livorno in
the year 1792, he showed his fully developed creative ingenuity in that part
of "sublime calculus" in which he was bound to find the subject for great
discoveries.

Called as Professor of Nautical Sciences in the College of Naval cadets in
Livorno in the year 1796 he published the Navigation Treatise of which were
printed three more editions more and more improved and detailed. This work
was and still is the only Italian textbook which is really suitable to
educate the practical pilot.

In the year 1798 was printed in Florence the work entitled "Calculus of
linear equations". In this oeuvre our author showed that he could
successfully compete with the most eminent geometers of Europe.

Postponing to a later discussion, as it will be suitable to do so while
talking about another book, the exposition of the many merits of this book,
we will limit ourselves here to say that while Laplace was calling falsely
not-integrable certain linear equations in which second order partial
differentials appear, while Paoli and Lacroix were investigating the same
subject and started to doubt about the statements of the mentioned French
geometer, Brunacci gave a method for integrating similar equations, being
able to generalize it to all differential orders. Paoli himself, by exposing
this method for a particular case in the third of the parts which form the
supplement to his Elements of Algebra, calls illustrious geometer that
scientist [i.e. Brunacci himself] who had been his student.

A voice finally was uttered from that place which had given birth to such
eminent scientists as Cavalieri, Frisi, Agnesi, and Oriani and Brunacci was
called to occupy the empty chair [of mathematics] in Pavia. He arrived there
in the year 1800 and although he arrived in a place where the Mathematical
Sciences were not ignored he met the greatest expectations and advanced the
fame [of that chair] to a yet unrivalled dignity.

Indeed it is not sufficient to be erudite in science to become its
professor, it is necessary to have the gift of the word, the capacity of
finding the right way to explain it. These gifts were given to him in the
highest, unrivalled level. Whoever heard him will admit that my expressions
although admired however are not enough to reveal the truth. The
mathematical teaching when coming from his lips was losing every difficulty
and bitterness, and developed with a peculiar charm and incantation [the
mathematical teaching] was at the same time education for the mind and
pleasure for the ears. It was then that the mathematical schools on the
banks of Ticino river reached the prestige which also nowadays is honoring
them. It was then that Vincenzo, having dedicated himself completely to his
science, started with all his forces to promote it.

The "Analisi Derivata" (Analysis of Derivatives) was printed in Pavia in the
year 1802. It is in this book that one can find one of the most sublime
concepts which was ever conceived by the human mind, that is the Principle
of Derivation. Because of it all the different parts of Mathematical
Sciences are tied and interconnected and it is opened an endless view which
allows\ us to consider as possible their infinite development. Soon he
conceived the challenging thought which lead him to re-write the whole body
of the doctrine of his science in many volumes, enriched by every novel
concept which had been formulated in the modern works. This endeavour may
have frightened everybody except him: he was also pushed by the advice of
that Sovereign Investigator of the stars who, being in Milan, wrote to
persuade him to start this oeuvre in the year 1800, believing that he was
the only one among the Italians who was capable to complete it successfully.

The oeuvre of the Course of Sublime Mathematics was printed in Florence in
four volumes in the years 1804, 1806, 1807, 1808. One would need a very long
time to expound, as it should be done, the merits of this book, but I will
want to shortly describe here its contents.

The first volume contains the Calculus of Finite Differences. This Calculus,
which was originated among the obscure calculations presented by Taylor,
which was developed in many Memoirs disseminated here and there in the
Proceedings of many Academies, for the first time was given scientific order
and method by the Florentine Geometer. He wrote it finding in his ingenuous
mind all that which was lacking in order to form a perfect theoretical
frame, and he infused in it all novels results which he had obtained in his
already mentioned works. It was his original contribution the integration of
linear equations of second order with variable coefficients, it was his
contribution a new formula for the integration of linear equation of all
orders with constant coefficients; it was his own the method to complete the
integrals to be replaced to the one proposed by D'Alembert, [method] which
he successfully introduced also in the differential calculus; but the idea
of the variable probability and the solution of the related problems, with
which he metaphorically could seize the wheel of the fortune and advanced in
the field where the genius of Lagragia had stopped, when in the Proceedings
of the Academy of Berlin (1775) he had given the solution of those problems
only in the case of constant probability. While citing the name of Lagrangia
I will not neglect to say that Brunacci was the first in Italy to see that
admirable light which the Theory of Analytical Functions can spread among
the mysterious smog which was obscuring the Infinitesimal Analysis. He
immediately conceived the idea of introducing it also among us: but oh! how
difficult was that endeavour! The Lagrangian notation, completely new,
produced a kind of revulsion: not all minds were firm enough to be able to
maintain -in the middle of a Revolution- their contact with the spirit of
the Calculus:\ He himself told me many times about the great obstacles which
he needed courageously to confront in pursuing his effort. He finally
managed to reach his aim, by reconciling the Lagrangian ideas with the
Leibnitz notation, together with the brackets introduced by Fontaine.

In this way are written the other three volumes of the said Course, where,
however, the author did not neglect to introduce with great skill whenever
possible the notation of the Geometer of Turin in order to make it familiar
to us.

We will only add that in the remaining part of his great oeuvre one can find
the rich results of Mathematical Analysis gathered from the most recent
Memoirs of the most celebrated Geometers and especially from the immense
body of works of the great Euler which he called his delight and from which
he admitted to have learnt that lucid order which makes his own works so
brilliant.

Oh! How many times I heard him talking about Euler with a great enthusiasm
and to urge me, and many of his other students, to study the work of the
only author who is suitable to educate a geometer! The great men, even when
are quoting the results of other authors are able to give to the subject
their own mark. This statement is true for that book where Brunacci infused
many of his ideas, not only those which we mentioned but also many others
which equally would merit to be mentioned, and in particular in those
various problems of every kind of applied mathematics and in the calculus of
variations which is reduced to the differential calculus and is there
exposed with a great detail, and finally in the mixed calculus, of which he
was the first to give the true principles and to expound in orderly way the
doctrine.

However a triumph which Brunacci obtained in front of all his rivals. The
Theory of the hydraulic water hammer, which seemed to be rebellious to the
lordship of Mathematical Analysis, and which was demanded with a golden
prize -without success- by the Academy of Berlin to the greatest geometers
of Europe in the year 1810 and then again in the year 1812 doubling the
prize, since 1810 was discovered by Brunacci who should have had received
the promised reward if an accident -which I do not want to recall here- had
not defrauded him of the deserved glory;\ this Theory was published in the
Treatise of the hydraulic water hammer of which were printed two editions;
in this Treatise said Theory, reduced to formulas and problems, is expounded
in the most efficacious way.

It is custom of the brave to prepare himself to the new victories and not to
be proud of the past ones. Therefore a new arena was chosen by our athlete
where he managed to defeat strong rivals. If he competes for discovering the
nature in hydraulic problems, Brunacci is awarded by the Societ\`{a}
Italiana: if he needs to reach the highest abstraction in order to find the
best metaphysics for the Calculus, Brunacci is awarded by the Accademia di
Padova. The Proceedings of the Illustrious Societ\`{a} Italiana carry the
name of Brunacci as author of many of the best Memoirs: too long would be to
cite all of them.

I will only mention that one on some particular solutions for the finite
difference equations, which our author treats in a way which is similar to
the one used by Lagrangia for differential equations, and where he
discovered some very elegant theorems valid for finite difference equations
which are not true for differential equations and the other one on shock
waves in fluids which embellishes the last Volume printed by said Societ\`{a}%
, Memoir where the analytical spirit really is triumphant.

Also the Istituto Nazionale Italiano was immediately honored in its first
Volume with a Memoir by Brunacci on the Theory of Maxima and High Minima;
subject which was remarkably advanced later in another Memoir. The Societ%
\`{a} Italiana and the Istituto oh! how greatly will grieve the loss of a
man who honored them with many and valuable works! Also the Academies of
Berlin, Munich, Turin and Lucca, and the others to which he belonged, will
perceive the great emptiness which is now left in them. I simply quickly
cite the textbook on the Elements of Algebra and Geometry written by our
author for the high school in few days, of which one has to praise the order
and the distribution of subjects and which was published in many editions.

I will mention as meriting great praise the Compendium of Sublime Calculus
which was issued in two Volumes in the year 1811, where it is gathered
everything which is sufficient to educate thoroughly a young geometer. In
writing it the author greatly improved and carefully modified many parts of
the complete course, and all added many new results and arguments.

It is not licit to neglect to indicate another subject in which -with
honored efforts- our professor distinguished himself. The Journal of
Physical Chemistry of Pavia was illustrated in many of his pages by his
erudite pen; I will content myself to indicate here three Memoirs where he
examines the doctrine of capillary attraction of Mister Laplace, comparing
it with that of Pessutti and where with his usual frankness, which is
originated by his being persuaded of how well-founded was his case, he
proves with his firm reasonings, whatever it is said\ by the French
geometers, some propositions which are of great praise for the mentioned
Italian geometer.

One could think that a man who wrote so much in his short life actually
should have been remained closed all the time alone in his office. On the
contrary: he not only was a great theoretician but also he was excellent in
all practical hydrometric and geodetic operations. He was Professor also in
these disciplines and with great dedication he worked heavily along the
banks of Ticino river in order to educate the best engineers.

He was a really skilled experimentalist and he often investigated natural
phenomena, getting favorable answers. I know very well how much interest
pushed him to these experimental activities, as is proven by the Hydrometry
Laboratory of the University which he founded and improved (sometimes at his
own expenses) with high quality instruments.

Also in these more practical activities his capabilities won him an
universal esteem, so that he was called everywhere sometimes on the river
banks in order to monitor their construction or for prevent their collapse
or sometimes on the navigation canals, among which the famous one in Pavia
was started under his direction which was confided to him by the past
government. The same government nominated him inspector of waters and
streets, inspector general of the public instruction and knight.

His character was strong in his resolutions, [it was] constant and resolute
in his sentiments, vigorous in the spirit, ready to well reason and ponder,
[it was] active and ready to engage in the [needed] efforts but above all he
was friendly and urbane: [his character] made him the center of social life
and the joy of friendship. Particularly with his students he was renouncing
to all the superiority of the "maestro" and assumed the attitude of the
father: I must avoid this memory, woe is me!, because it too strongly makes
tears to come to my eyes. Those who need the evidence of my last statement
has simply to see his how his students wanted to honor that great man and to
manifest their sorrow: they carried on their shoulders his mortal remains,
they decorated in an extraordinary way his funeral parlour and now are
praising the departed's merits with their tears and their silent grief which
are more eloquent than all spoken lamentation.

Everybody who is now promising to contribute to exact sciences in Lombardy
is a student of Brunacci, and indeed among his disciples there are those
who, as their mentor himself often said, is now an eagle who can fly with
his own wings. Such [an eagle] is the Professor in Bologna, author of the
essay on Poligonometry, such is the other one who is the author of the
Treatise on the Contours of the Shadows and whose noteworthy voice is
entitled to succeed to that of his Maestro on the banks of Ticino, such is a
third disciple who has already shown that Italy can hope to have soon a
Geometer who will emulate the great genius who wrote the theory of celestial
bodies.

What a great misfortune was to see the departure of a man in the age of his
maturity who already had greatly contributed to science and who was bound to
contribute even more copiously to it! I know very well, as I had many times
the privilege of his confidences about the subject of his studies, how many
precious works can be found in his manuscripts. Among them, some excellent
documents which he wanted to gather to form a commentary to the Analytical
Mechanics, many very beautiful discourses read on occasion of the defense of
theses, some sequels of Memoirs containing the description and the
calculation of many machines inspired by the Hydraulic Architecture authored
by Belidor, gathering which he intended to complete an oeuvre which would
have been of great utility.

May these last achievements of such an inventive and ingenious Geometer be
delivered up to a capable and educated scholar, who could enlighten them as
they deserve, for the advancement of SCIENCES, for the glory of the AUTHOR\
and for the prestige of ITALY

Milan, 18 June 1818

\subsection{Appendix C. Peridynamics:\ A new/old model for deformable bodies}

The celebrated and fundamental textbook by Lagrange \cite{lagrange1788} is,
with few and biased exceptions, generally regarded as a milestone in
Mechanical Sciences and unanimously as novel in its content and style of
presentation.

Indeed Lagrange himself, differently from what was done by his epigones,
puts his work in the correct perspective, by giving the due credit to all
his predecessors. Indeed the M\'{e}canique Analytique starts with an
interesting historical introduction, which can be considered the initiation
of the modern history of mechanics. Unfortunately also this aspect of the
Lagrangian lesson is not very often followed in modern science.

A very new Continuum Mechanical Theory has been recently announced and
developed: Peridynamics. Actually the ideas underlying Peridynamics are very
interesting and most likely they deserve the full attention of experts in
continuum, fracture and damage mechanics.

Indeed starting from a balance law of the form (\ref{Peridynamics}) for
instance in \cite{dipaolaetal2010}, \cite{dipaolaetal2010a} and \cite%
{silling2000} (but many other similar treatments are available in the
literature) one finds a formulation of Continuum Mechanics which relaxes the
standard one transmitted by the apologists of the Cauchy format and seems
suitable (see the few comments below) to describe many and interesting
phenomena e.g. in crack formation and growth.

However even those scientists whose mother language is Italian actually seem
unaware of the contribution due to Gabrio Piola in this field: this loss of
memory and this lack of credit to the major sources of our knowledge, even
in those cases in which their value is still topical, is very dangerous, as
proven in detail by the analysis developed in Russo \cite{russo2013}, \cite%
{russo2003}.

Unfortunately this tendency towards a mindless "modernism" seems to become
more and more aggravated.

In \cite{silling2000} the analysis started by Piola is continued, seemingly
as if the author, Silling, were one of his closer pupils: arguments are very
similar and also a variational formulation of the presented theories is
found and discussed. In \cite{lehoucqsilling2008} and in \cite%
{sillinglehoucq2008} it is stated in the abstracts that:

``The peridynamic model is a framework for
continuum mechanics based on the idea that pairs of particles exert forces
on each other across a finite distance. The equation of motion in the
peridynamic model is an integro- differential equation. In this paper, a
notion of a peridynamic stress tensor derived from nonlocal interactions is
defined.''

``The peridynamic model of solid mechanics is a
nonlocal theory containing a length scale. It is based on direct
interactions between points in a continuum separated from each other by a
finite distance. The maximum interaction distance provides a length scale
for the material model. This paper addresses the question of whether the
peridynamic model for an elastic material reproduces the classical local
model as this length scale goes to zero. We show that if the motion,
constitutive model, and any nonhomogeneities are sufficiently smooth, then
the peridynamic stress tensor converges in this limit to a Piola-Kirchhoff
stress tensor that is a function only of the local deformation gradient
tensor, as in the classical theory. This limiting Piola-Kirchhoff stress
tensor field is differentiable, and its divergence represents the force
density due to internal forces.''

The reader is invited to compare these statements with those which can be
found in the previous Appendix A.

It is very interesting to see how fruitful can be the ideas formulated 167
years ago by Piola. It is enough to read the abstract of \cite%
{askarietal2008}

``The paper presents an overview of peridynamics,
a continuum theory that employs a nonlocal model of force interaction.
Specifically, the stress/strain relationship of classical elasticity is
replaced by an integral operator that sums internal forces separated by a
finite distance. This integral operator is not a function of the deformation
gradient, allowing for a more general notion of deformation than in
classical elasticity that is well aligned with the kinematic assumptions of
molecular dynamics. Peridynamics' effectiveness has been demonstrated in
several applications, including fracture and failure of composites,
nanofiber networks, and polycrystal fracture. These suggest that
peridynamics is a viable multiscale material model for length scales ranging
from molecular dynamics to those of classical elasticity.''

Or also the abstract of the paper by Parks et al. \cite{parksetal2008}.

``Peridynamics (PD) is a continuum theory that
employs a nonlocal model to describe material properties. In this context,
nonlocal means that continuum points separated by a finite distance may
exert force upon each other. A meshless method results when PD is
discretized with material behavior approximated as a collection of
interacting particles. This paper describes how PD can be implemented within
a molecular dynamics (MD) framework, and provides details of an efficient
implementation. This adds a computational mechanics capability to an MD code
enabling simulations at mesoscopic or even macroscopic length and time scales
''

It is remarkable how strictly related are non-local continuum theories with
the discrete theories of particles bound to the nodes of a lattice. How deep
was the insight of Piola can be understood by looking at the literature
about the subject which includes for instance \cite{askarietal2008}, \cite%
{demmiesilling2007}, \cite{dipaolaetal2010}, \cite{dipaolaetal2010a}, \cite%
{duetal2013}, \cite{emmrichetal2013}, \cite{lehoucqsilling2008}, \cite%
{selesonetal2013}, \cite{silling2000}, \cite{sillingetal2007}, \cite%
{sillinglehoucq2008}.

\subsection{Appendix D. On an expression for $\protect\nabla F$ deduced in
Piola (1845-6) on pages 158-159}

In this appendix we deduce, by means of the Levi-Civita tensor calculus, the
expression for the second gradient of placement that is needed to transform
eqn. (\ref{FirstDerivative}) into eqn. (\ref{THIRDDERIVATIVE}) and that is
obtained by \cite{piola1845-6}, see the Appendix A for the appropriate
translation. The original calculations are rather lengthy and cumbersome: it
is however the opinion of the authors that Piola had caught their
"tensorial" or at least their algebraic structure. Indeed the notation he
used made rather easy the identification of the tensorial objects involved.

We start from the following identification between modern and Piola's
notation

\begin{equation}
\begin{array}{ccc}
F_{\alpha }^{i}\leftrightarrows \left( 
\begin{array}{ccc}
\frac{dx}{da} & \frac{dy}{da} & \frac{dz}{da} \\ 
\frac{dx}{db} & \frac{dy}{db} & \frac{dz}{db} \\ 
\frac{dx}{dc} & \frac{dy}{dc} & \frac{dz}{dc}%
\end{array}%
\right) & \det F\leftrightarrows H & \left( \det F\right) \left(
F^{-1}\right) _{\text{ \hspace{0in} \hspace{0in}}j}^{\beta }\leftrightarrows
\left( 
\begin{array}{ccc}
l_{1} & m_{1} & n_{1} \\ 
l_{2} & m_{2} & n_{2} \\ 
l_{3} & m_{3} & n_{3}%
\end{array}%
\right)%
\end{array}
\tag{N7}  \label{N10bisbisbis}
\end{equation}%
so that we can state that equation (28) on page 26 of \cite{piola1845-6} is
equivalent to the following one 
\begin{equation*}
\left( \det F\right) F^{-1}F=\left( \det F\right) I
\end{equation*}%
where $I$ is the identity matrix. Moreover the equation (6) on page 57 is
equivalent to the following one 
\begin{equation}
\begin{array}{cc}
\left( 
\begin{array}{ccc}
t_{1} & t_{4} & t_{5} \\ 
t_{4} & t_{2} & t_{6} \\ 
t_{5} & t_{6} & t_{3}%
\end{array}%
\right) \leftrightarrows C=F^{T}F & C_{\alpha \beta }=F_{\alpha
}^{i}F_{i\beta }%
\end{array}
\tag{N8}  \label{N10bisbis}
\end{equation}

The equation on pages 158-159 in \cite{piola1845-6} is written in tensorial
form as follows:%
\begin{equation}
\begin{array}{cc}
\frac{\partial ^{2}\chi ^{i}}{\partial X^{\alpha }\partial X^{\beta }}%
:=D_{\alpha \beta \eta }\left( F^{-1}\right) ^{i\eta }, & F_{i\eta }\frac{%
\partial ^{2}\chi ^{i}}{\partial X^{\alpha }\partial X^{\beta }}=\frac{%
\partial F_{\alpha }^{i}}{\partial X^{\beta }}F_{i\eta }=D_{\alpha \beta
\eta }%
\end{array}
\tag{N9}  \label{N10bis}
\end{equation}

Now by recalling that 
\begin{equation*}
\frac{\partial F_{i\beta }}{\partial X^{\gamma }}=\frac{\partial F_{i\gamma }%
}{\partial X^{\beta }}
\end{equation*}%
we have the symmetry of $D$ with respect to the first two indeces,%
\begin{equation}
D_{\alpha \beta \eta }=D_{\beta \alpha \eta }  \tag{N10}  \label{N11bis}
\end{equation}%
and, because of such expression, we can perform the following simple
calculations (usual symmetrization, $A_{(ab)}=A_{ab}+A_{ba}$, and
skew-symmetrization, $A_{[ab]}=A_{ab}-A_{ba}$, conventional symbols are used)

\begin{equation*}
\frac{\partial C_{\alpha\beta}}{\partial X^{\gamma}}=\frac{\partial
F_{\alpha }^{i}}{\partial X^{\gamma}}F_{i\beta}+F_{\alpha}^{i}\frac{\partial
F_{i\beta}}{\partial X^{\gamma}}
\end{equation*}

\begin{align*}
2\frac{\partial C_{\alpha\lbrack\beta}}{\partial X^{\gamma]}} & =\frac{%
\partial C_{\alpha\beta}}{\partial X^{\gamma}}-\frac{\partial
C_{\alpha\gamma}}{\partial X^{\beta}}= \\
& =\frac{\partial F_{\alpha}^{i}}{\partial X^{\gamma}}F_{i\beta}+F_{\alpha
}^{i}\frac{\partial F_{i\beta}}{\partial X^{\gamma}}-\frac{\partial
F_{\alpha }^{i}}{\partial X^{\beta}}F_{i\gamma}-F_{\alpha}^{i}\frac{\partial
F_{i\gamma }}{\partial X^{\beta}}= \\
& =\frac{\partial F_{\alpha}^{i}}{\partial X^{\gamma}}F_{i\beta}-\frac{%
\partial F_{\alpha}^{i}}{\partial X^{\beta}}F_{i\gamma}=D_{\alpha
\gamma\beta}-D_{\alpha\beta\gamma}=2D_{\alpha\lbrack\gamma\beta]}
\end{align*}

\begin{align*}
\frac{\partial C_{\alpha\beta}}{\partial X^{\gamma}} & =\frac{\partial\left(
F_{\alpha}^{i}F_{i\beta}\right) }{\partial X^{\gamma}}=\frac{\partial
F_{\alpha}^{i}}{\partial X^{\gamma}}F_{i\beta}+F_{\alpha}^{i}\frac{\partial
F_{i\beta}}{\partial X^{\gamma}}= \\
& =D_{\gamma\alpha\beta}+D_{\gamma\beta\alpha}=2D_{\gamma(\alpha\beta)}
\end{align*}
By decomposing $D$ into its skew and symmetric parts (with respect to the
second and third index, see also (\ref{N11bis})) one gets 
\begin{equation}
D_{\gamma\alpha\beta}=D_{\gamma(\alpha\beta)}+D_{\gamma\lbrack\alpha\beta ]}=%
\frac{1}{2}\frac{\partial C_{\alpha\beta}}{\partial X^{\gamma}}+\frac{%
\partial C_{\gamma\lbrack\beta}}{\partial X^{\alpha]}}=\frac{1}{2}\left( 
\frac{\partial C_{\alpha\beta}}{\partial X^{\gamma}}+\frac{\partial
C_{\gamma\beta}}{\partial X^{\alpha}}-\frac{\partial C_{\gamma\alpha}}{%
\partial X^{\beta}}\right)  \tag{N11}  \label{RepresentationD}
\end{equation}

The third order tensor $D_{\gamma\alpha\beta}$ which we have introduced
allows us to reproduce in the compact form (\ref{N10bis}) the formula which
occupies nearly two pages of Piola's work. Moreover we have obtained the
formula \ref{RepresentationD} with an easy calculation process which is much
less involved than the one first conceived by Piola.

From (\ref{RepresentationD}) we have

\begin{equation}
F_{i\beta }\frac{\partial ^{2}\chi ^{i}}{\partial X^{\alpha }\partial
X^{\gamma }}=\frac{1}{2}\left( \frac{\partial C_{\alpha \beta }}{\partial
X^{\gamma }}+\frac{\partial C_{\gamma \beta }}{\partial X^{\alpha }}-\frac{%
\partial C_{\gamma \alpha }}{\partial X^{\beta }}\right)  \tag{N12}
\label{objectivityC}
\end{equation}%
which is equivalent to 
\begin{equation}
\frac{\partial ^{2}\chi ^{j}}{\partial X^{\alpha }\partial X^{\gamma }}=%
\frac{1}{2}\left( F^{-1}\right) ^{j\beta }\left( \frac{\partial C_{\alpha
\beta }}{\partial X^{\gamma }}+\frac{\partial C_{\gamma \beta }}{\partial
X^{\alpha }}-\frac{\partial C_{\gamma \alpha }}{\partial X^{\beta }}\right) 
\tag{N13}  \label{gradFPiola}
\end{equation}%
To compare the two formalisms let us state the identification of the
left-hand side of one line, i.e. of the 11th one divided by $2H$ of the
formula appearing on page 158 in \cite{piola1845-6}, i.e.,

\begin{equation*}
\frac{\partial^{2}\chi^{2}}{\partial X^{1}\partial X^{2}}\leftrightarrows 
\text{ }\frac{d^{2}y}{dadb}.
\end{equation*}
Thus, from (\ref{RepresentationD}) with $\alpha=1$, $j=\gamma=2$, by
recalling the symmetry of the tensor $C$ and the identifications (\ref%
{N10bisbisbis}) and (\ref{N10bisbis}),

\begin{gather*}
\frac{\partial ^{2}\chi ^{2}}{\partial X^{1}\partial X^{2}}=\frac{1}{2}%
\left( F^{-1}\right) ^{2\beta }\left( \frac{\partial C_{1\beta }}{\partial
X^{2}}+\frac{\partial C_{2\beta }}{\partial X^{1}}-\frac{\partial C_{21}}{%
\partial X^{\beta }}\right) = \\
=\frac{1}{2}\left( F^{-1}\right) ^{21}\left( \frac{\partial C_{11}}{\partial
X^{2}}+\frac{\partial C_{21}}{\partial X^{1}}-\frac{\partial C_{21}}{%
\partial X^{1}}\right) + \\
+\frac{1}{2}\left( F^{-1}\right) ^{22}\left( \frac{\partial C_{12}}{\partial
X^{2}}+\frac{\partial C_{22}}{\partial X^{1}}-\frac{\partial C_{21}}{%
\partial X^{2}}\right) + \\
+\frac{1}{2}\left( F^{-1}\right) ^{23}\left( \frac{\partial C_{13}}{\partial
X^{2}}+\frac{\partial C_{23}}{\partial X^{1}}-\frac{\partial C_{21}}{%
\partial X^{3}}\right) = \\
=\frac{1}{2}\left( F^{-1}\right) ^{21}\frac{\partial C_{11}}{\partial X^{2}}+%
\frac{1}{2}\left( F^{-1}\right) ^{22}\frac{\partial C_{22}}{\partial X^{1}}+
\\
+\frac{1}{2}\left( F^{-1}\right) ^{23}\left( \frac{\partial C_{13}}{\partial
X^{2}}+\frac{\partial C_{23}}{\partial X^{1}}-\frac{\partial C_{21}}{%
\partial X^{3}}\right) \leftrightarrows \\
\leftrightarrows \frac{l_{2}}{2H}\frac{dt_{1}}{db}+\frac{m_{2}}{2H}\frac{%
dt_{2}}{da}+\frac{n_{2}}{2H}\left( \frac{dt_{5}}{db}+\frac{dt_{6}}{da}-\frac{%
dt_{4}}{dc}\right)
\end{gather*}%
which is, multiplying by $2H$ both members, the 11th equality on page 158 in 
\cite{piola1845-6}%
\begin{equation*}
2H\text{ }\frac{d^{2}y}{dadb}=l_{2}\frac{dt_{1}}{db}+m_{2}\frac{dt_{2}}{da}%
+n_{2}\left( \frac{dt_{5}}{db}+\frac{dt_{6}}{da}-\frac{dt_{4}}{dc}\right)
\end{equation*}

Piola continued the calculations by considering the third order derivatives.
However the obtained expressions are too long for being reproduced in
printed form. So he states:

``The trinomials with third order derivatives are
of three kinds: there are those constituted by derivatives of first and
third order, and one can count 30 of them: there are those constituted by
derivatives of second and third order, and they are 60 in number: and there
are those which contain only third order derivatives and they are 55 in
number. I am not writing them, as everybody who is given the needed patience
can easily calculate them by himself, as it can be also done for those
trinomials containing derivatives of higher order.''

As we can use Levi-Civita tensor calculus it is easier for us to find the
needed patience, at least for calculating the trinomials constituted by
derivatives of first and third order. Indeed from

\begin{equation}
F_{i\gamma}\frac{\partial^{2}\chi^{i}}{\partial X^{\alpha}\partial X^{\beta}}%
=\frac{1}{2}\left( \frac{\partial C_{\alpha\gamma}}{\partial X^{\beta}}+%
\frac{\partial C_{\beta\gamma}}{\partial X^{\alpha}}-\frac{\partial
C_{\beta\alpha}}{\partial X^{\gamma}}\right)  \tag{N14}
\label{RiemannIdentity}
\end{equation}
by differentiating the (\ref{RiemannIdentity}) we get%
\begin{equation*}
\frac{\partial}{\partial X^{\eta}}\left( F_{i\gamma}\frac{%
\partial^{2}\chi^{i}}{\partial X^{\alpha}\partial X^{\beta}}\right) =\frac{%
\partial }{\partial X^{\eta}}\left( \frac{1}{2}\left( \frac{\partial
C_{\alpha\gamma }}{\partial X^{\beta}}+\frac{\partial C_{\beta\gamma}}{%
\partial X^{\alpha}}-\frac{\partial C_{\beta\alpha}}{\partial X^{\gamma}}%
\right) \right)
\end{equation*}
and rearranging the terms,

\begin{equation}
F_{i\gamma}\frac{\partial^{3}\chi^{i}}{\partial X^{\alpha}\partial X^{\beta
}\partial X^{\eta}}=-\frac{\partial^{2}\chi_{i}}{\partial X^{\gamma}\partial
X^{\eta}}\frac{\partial^{2}\chi^{i}}{\partial X^{\alpha}\partial X^{\beta}}+%
\frac{1}{2}\left( \frac{\partial^{2}C_{\alpha\gamma}}{\partial X^{\eta
}\partial X^{\beta}}+\frac{\partial^{2}C_{\beta\gamma}}{\partial X^{\eta
}\partial X^{\alpha}}-\frac{\partial^{2}C_{\beta\alpha}}{\partial X^{\eta
}\partial X^{\gamma}}\right) .  \tag{N15}  \label{N14bis}
\end{equation}

By replacing the following equality due to (\ref{N10bis})%
\begin{equation}
\frac{\partial^{2}\chi^{i}}{\partial X^{\alpha}\partial X^{\beta}}=\frac{1}{2%
}\left( F^{-1}\right) ^{i\delta}\left( \frac{\partial C_{\alpha\delta}}{%
\partial X^{\beta}}+\frac{\partial C_{\beta\delta}}{\partial X^{\alpha}}-%
\frac{\partial C_{\beta\alpha}}{\partial X^{\delta}}\right)  \tag{N16}
\label{SecondDerivative}
\end{equation}
in the identity (\ref{N14bis}) one gets%
\begin{align*}
F_{i\gamma}\frac{\partial^{3}\chi^{i}}{\partial X^{\alpha}\partial X^{\beta
}\partial X^{\eta}} & =-\frac{1}{2}\left( F^{-1}\right) _{i}^{\nu}\left( 
\frac{\partial C_{\gamma\nu}}{\partial X^{\eta}}+\frac{\partial C_{\eta\nu}}{%
\partial X^{\gamma}}-\frac{\partial C_{\eta\gamma}}{\partial X^{\nu}}\right) 
\frac{1}{2}\left( F^{-1}\right) ^{i\delta}\left( \frac{\partial
C_{\alpha\delta}}{\partial X^{\beta}}+\frac{\partial C_{\beta\delta}}{%
\partial X^{\alpha}}-\frac{\partial C_{\beta\alpha}}{\partial X^{\delta}}%
\right) \\
& +\frac{1}{2}\left( \frac{\partial^{2}C_{\alpha\gamma}}{\partial X^{\eta
}\partial X^{\beta}}+\frac{\partial^{2}C_{\beta\gamma}}{\partial X^{\eta
}\partial X^{\alpha}}-\frac{\partial^{2}C_{\beta\alpha}}{\partial X^{\eta
}\partial X^{\gamma}}\right)
\end{align*}
which can easily be rewritten in the form

\begin{align*}
F_{i\gamma}\frac{\partial^{3}\chi^{i}}{\partial X^{\alpha}\partial X^{\beta
}\partial X^{\eta}} & =-\frac{1}{4}\left( C^{-1}\right) ^{\nu\delta }\left( 
\frac{\partial C_{\gamma\nu}}{\partial X^{\eta}}+\frac{\partial C_{\eta\nu}}{%
\partial X^{\gamma}}-\frac{\partial C_{\eta\gamma}}{\partial X^{\nu}}\right)
\left( \frac{\partial C_{\alpha\delta}}{\partial X^{\beta}}+\frac{\partial
C_{\beta\delta}}{\partial X^{\alpha}}-\frac{\partial C_{\beta\alpha}}{%
\partial X^{\delta}}\right) \\
& +\frac{1}{2}\left( \frac{\partial^{2}C_{\alpha\gamma}}{\partial X^{\eta
}\partial X^{\beta}}+\frac{\partial^{2}C_{\beta\gamma}}{\partial X^{\eta
}\partial X^{\alpha}}-\frac{\partial^{2}C_{\beta\alpha}}{\partial X^{\eta
}\partial X^{\gamma}}\right)
\end{align*}
which has the structure sought after by Piola.

With easy calculation, from the last equation we get%
\begin{equation}
\frac{1}{2}\frac{\partial^{3}\rho^{2}(\bar{X},X)}{\partial\bar{X}_{\alpha
}\partial\bar{X}_{\beta}\partial\bar{X}_{\gamma}}=\frac{1}{2}\left( \frac{%
\partial C_{\alpha\gamma}}{\partial X^{\beta}}+\frac{\partial C_{\beta\gamma}%
}{\partial X^{\alpha}}+\frac{\partial C_{\beta\alpha}}{\partial X^{\gamma}}%
\right) +\left( \sum\limits_{i=1}^{3}\left( \chi_{i}(\bar {X}%
)-\chi_{i}(X)\right) \frac{\partial^{3}\chi_{i}(\bar{X})}{\partial\bar {X}%
_{\alpha}\partial\bar{X}_{\beta}\partial\bar{X}_{\gamma}}\right)  \tag{N17}
\label{ThirdDerivativeRHO}
\end{equation}
In the following appendix E an induction argument will be presented which
allows us to prove the conjecture put forward by Piola at the beginning of
his num.74 pag.156.

\subsection{Appendix E. "After these calculations I abandoned myself to the
analogy"}

We are not so enthusiastic about the work of Piola to the extent that we
cannot see clearly the limits of his mathematical proofs.

Indeed the important property which he discusses in the num.74 is obtained
by means of a proof "by analogy" which is not considered acceptable
nowadays. Although there are examples of mathematical induction which are
very ancient (see the discussion in \cite{russo2003} and references therein)
only after Boole and Dedekind it became a universally known and (nearly
universally) accepted method.

Actually Piola states here that, because of objectivity, the expression of
Virtual Work must depend only on deformation measure $C_{\gamma \beta }$ and
its derivatives. However, as we have already pointed out, his proof is
based, for higher derivatives, on an argument which the majority of
contemporary mathematicians would consider no more than a (maybe
well-grounded) conjecture. Indeed at the beginning of page 157 of \cite%
{piola1845-6} one reads

"\textbf{after these calculations I abandoned myself to the analogy: and
this will be sooner or later unavoidable, because our series is infinite and
it will be impossible to check all its terms."}

We reproduce here an inductive argument which indeed follows the original
spirit of Piola.

Let us start by proving that:

\begin{lemma}
\textbf{Representation of placement higher order derivatives}. For every $n$
there exist a family of (polynomial) functions $M_{\gamma\alpha_{1}...%
\alpha_{n}}$ of the tensor variables $C,\nabla C,..\nabla^{n-1}C$ such that%
\begin{equation}
\left( \frac{\partial\chi_{i}(\bar{X})}{\partial\bar{X}^{\gamma}}\frac{%
\partial^{n}\chi^{i}(\bar{X})}{\partial\bar{X}^{\alpha_{1}}...\partial\bar{X}%
^{\alpha_{n}}}\right) =M_{\gamma\alpha_{1}...\alpha_{n}}\left(
C,..,\nabla^{n-1}C\right)  \tag{N18}  \label{HYPREC1}
\end{equation}
\end{lemma}

As in the Appendix D we have proven such a lemma for $n=2$ that is,%
\begin{equation}
F_{i\gamma}\frac{\partial^{2}\chi^{i}}{\partial X^{\alpha}\partial X^{\beta}}%
=\frac{1}{2}\left( \frac{\partial C_{\alpha\gamma}}{\partial X^{\beta}}+%
\frac{\partial C_{\beta\gamma}}{\partial X^{\alpha}}-\frac{\partial
C_{\beta\alpha}}{\partial X^{\gamma}}\right) .  \tag{N19}
\label{HYPREC1-N=2}
\end{equation}
In order to prove (\ref{HYPREC1}) for every $n$ it is sufficient to prove
that if it is valid for all $N\leq n$ then it is valid also for $N=n+1.$

Let us start by remarking that (\ref{HYPREC1}) implies that

\begin{equation}
\frac{\partial^{n}\chi^{i}(\bar{X})}{\partial\bar{X}^{\alpha_{1}}...\partial%
\bar{X}^{\alpha_{n}}}=\left( F^{-1}\right) ^{i\eta}M_{\eta
\alpha_{1}...\alpha_{n}}\left( C,..,\nabla^{n-1}C\right)  \tag{N20}
\label{HYPREC1INVERSE}
\end{equation}

Let us then differentiate (\ref{HYPREC1}) assumed valid for $N=n$ to get 
\begin{equation*}
\frac{\partial}{\partial\bar{X}^{\alpha_{n+1}}}\left( \frac{\partial\chi
_{i}(\bar{X})}{\partial\bar{X}^{\gamma}}\frac{\partial^{n}\chi^{i}(\bar{X})}{%
\partial\bar{X}^{\alpha_{1}}...\partial\bar{X}^{\alpha_{n}}}\right) =\frac{%
\partial}{\partial\bar{X}^{\alpha_{n+1}}}\left( M_{\gamma\alpha
_{1}...\alpha_{n}}\left( C,..,\nabla^{n-1}C\right) \right)
\end{equation*}
which implies 
\begin{align*}
\frac{\partial\chi_{i}(\bar{X})}{\partial\bar{X}^{\gamma}}\frac{\partial
^{n+1}\chi^{i}(\bar{X})}{\partial\bar{X}^{\alpha_{1}}...\partial\bar {X}%
^{\alpha_{n}}\partial\bar{X}^{\alpha_{n+1}}} & =\frac{\partial}{\partial\bar{%
X}^{\alpha_{n+1}}}\left( M_{\gamma\alpha_{1}...\alpha_{n}}\left(
C,..,\nabla^{n-1}C\right) \right) \\
& -\frac{\partial^{2}\chi_{i}(\bar{X})}{\partial\bar{X}^{\gamma}\partial 
\bar{X}^{\alpha_{n+1}}}\frac{\partial^{n}\chi^{i}(\bar{X})}{\partial\bar {X}%
^{\alpha_{1}}...\partial\bar{X}^{\alpha_{n}}}
\end{align*}

Now by replacing equation (\ref{HYPREC1INVERSE}) two times (for $n=2$ and
for $N=n)$ we get

\begin{gather*}
\frac{\partial\chi_{i}(\bar{X})}{\partial\bar{X}^{\gamma}}\frac{\partial
^{n+1}\chi^{i}(\bar{X})}{\partial\bar{X}^{\alpha_{1}}...\partial\bar {X}%
^{\alpha_{n}}\partial\bar{X}^{\alpha_{n+1}}}=\frac{\partial}{\partial \bar{X}%
^{\alpha_{n+1}}}\left( M_{\gamma\alpha_{1}...\alpha_{n}}\left(
C,..,\nabla^{n-1}C\right) \right) \\
-\left( \left( F^{-1}\right) _{i}^{\eta}M_{\eta\alpha_{n+1}}\left(
C,..,\nabla^{n-1}C\right) \right) \left( \left( F^{-1}\right) ^{i\beta
}M_{\beta\alpha_{1}...\alpha_{n}}\left( C,..,\nabla^{n-1}C\right) \right) =
\\
=M_{\gamma\alpha_{1}...\alpha_{n}\alpha_{n+1}}\left( C,..,\nabla
^{n-1}C\right) ,
\end{gather*}
where we have introduced the definition 
\begin{align*}
M_{\gamma\alpha_{1}...\alpha_{n}\alpha_{n+1}}\left( C,..,\nabla
^{n-1}C\right) & :=\frac{\partial}{\partial\bar{X}^{\alpha_{n+1}}}\left(
M_{\gamma\alpha_{1}...\alpha_{n}}\left( C,..,\nabla^{n-1}C\right) \right) \\
& -\left( C^{-1}\right) ^{\beta\eta}M_{\eta\alpha_{n+1}}\left(
C,..,\nabla^{n-1}C\right) M_{\beta\alpha_{1}...\alpha_{n}}\left(
C,..,\nabla^{n-1}C\right)
\end{align*}

The proof by induction of the\textbf{\ }lemma is thus complete. To prove
that also the generic $n-th$ order derivative of $\rho ^{2}$ can be
expressed, when $\bar{X}=X$ , in terms of the $\left( n-2\right) -th$ order
derivatives of $C_{\gamma \beta }$ we can use again a simple recursion
argument based on the previous lemma.

Indeed the following other lemma is true:

\begin{lemma}
\textbf{Representation of the derivatives of the distance function }$\rho.$
For every $n$ there exist a family of (polynomial) functions $L_{\alpha
_{1}....\alpha_{n}}$ of the variables $C,..,\nabla^{n-2}C$ such that 
\begin{equation}
\frac{\partial^{n}\rho^{2}(\bar{X},X)}{\partial\bar{X}^{\alpha_{1}}...%
\partial\bar{X}^{\alpha_{n}}}=L_{\alpha_{1}....\alpha_{n}}\left(
C,..,\nabla^{n-2}C\right) +\left( \left( \chi_{i}(\bar{X})-\chi
_{i}(X)\right) \frac{\partial^{n}\chi^{i}(\bar{X})}{\partial\bar{X}%
^{\alpha_{1}}...\partial\bar{X}^{\alpha_{n}}}\right)  \tag{N21}
\label{H1RIC}
\end{equation}
\end{lemma}

To prove the Lemma we assume by inductive hypothesis that it is true for $%
N=n $ and prove that it is true for $N=n+1.$ As we have proven formula (\ref%
{ThirdDerivativeRHO}), that is the previous lemma for $n=3$, then the Lemma
follows by the Mathematical Induction Principle.

Therefore by differentiating equation (\ref{H1RIC}) one gets 
\begin{align*}
\frac{\partial^{n+1}\rho^{2}(\bar{X},X)}{\partial\bar{X}^{\alpha_{1}}...%
\partial\bar{X}^{\alpha_{n+1}}} & =\frac{\partial}{\partial\bar {X}%
^{\alpha_{n+1}}}\left( L_{\alpha_{1}....\alpha_{n}}\left( C,..,\nabla
^{n-2}C\right) \right) +\left( \frac{\partial\chi_{i}(\bar{X})}{\partial\bar{%
X}^{\alpha_{n+1}}}\frac{\partial^{n}\chi^{i}(\bar{X})}{\partial\bar{X}%
^{\alpha_{1}}...\partial\bar{X}^{\alpha_{n}}}\right) \\
& +\left( \sum\limits_{i=1}^{3}\left( \chi_{i}(\bar{X})-\chi_{i}(X)\right) 
\frac{\partial^{n+1}\chi_{i}(\bar{X})}{\partial\bar{X}^{\alpha_{1}}...%
\partial\bar{X}^{\alpha_{n}}\partial\bar{X}^{\alpha_{n+1}}}\right)
\end{align*}
which by replacing equation (\ref{HYPREC1}) becomes 
\begin{align*}
\frac{\partial^{n+1}\rho^{2}(\bar{X},X)}{\partial\bar{X}^{\alpha_{1}}...%
\partial\bar{X}^{\alpha_{n+1}}} & =\frac{\partial}{\partial\bar {X}%
^{\alpha_{n+1}}}\left( L_{\alpha_{1}....\alpha_{n}}\left( C,..,\nabla
^{n-2}C\right) \right) +M_{\alpha_{n+1}\alpha_{1}...\alpha_{n}}\left(
C,..,\nabla^{n-1}C\right) \\
& +\left( \sum\limits_{i=1}^{3}\left( \chi_{i}(\bar{X})-\chi_{i}(X)\right) 
\frac{\partial^{n+1}\chi_{i}(\bar{X})}{\partial\bar{X}^{\alpha_{1}}...%
\partial\bar{X}^{\alpha_{n}}\partial\bar{X}^{\alpha_{n+1}}}\right)
\end{align*}
which proves the Lemma once one has introduced the following recursive
definition 
\begin{equation*}
L_{\alpha_{1}....\alpha_{n}\alpha_{n+1}}\left( C,..,\nabla^{n-2}C\right) :=%
\frac{\partial}{\partial\bar{X}^{\alpha_{n+1}}}\left( L_{\alpha
_{1}....\alpha_{n}}\left( C,..,\nabla^{n-2}C\right) \right) +M_{\alpha
_{n+1}\alpha_{1}...\alpha_{n}}\left( C,..,\nabla^{n-1}C\right) .
\end{equation*}

\subsection{Appendix F An Italian Mathematical Genealogy}

In \cite{russo2013} it is discussed a most likely loss of knowledge occurred
at the end of the Punic wars. More generally all the processes of erasure
and removal of previosly well-established scientific knowledge are related
to the simultaneous occurrence of two circumstances: the loss of continuity
in the chain between Maestro and student in the academic institutions and
the loss of the awareness -in the whole society- of the strict connection
which exists between science (in all its most abstract expressions,
including mathematics) and technology (see also \cite{russo2003}).

The final result of the simultaneous occurrence of these two circumstances
is that the societies in which they occur do not invest resources in the
storage and transmission of theoretical knowledge and that, as a
consequence, it is broken the contact between maestro and pupil, established
when a living scientist teaches to his students the content of the most
difficult and important textbooks. As a final result, in those societies, at
first the theoretical knowledge, and subsequently after a more or less long
time period, also the technological capabilities are lost.

We want to underline in this appendix that there is a direct genealogy
starting from Gabrio Piola and leading to the founders of absolute tensor
calculus. The Italian school of the XIX Century was started under the
momentum impressed by the Napoleonic reforms of the political organization
of the Italian Nation: in this context the reader should see the Eulogy
in Appendix B where Piola, talking about the textbook of Mathematical
Analysis written by his Maestro Vincenzo Brunacci writes

``was also pushed by the advice of that Sovereign
who was an investigator of the stars who, being in Milan, wrote to persuade
him to start this oeuvre in the year 1800, believing that he was the only
one among the Italians who was capable to complete it successfully.''.

Gabrio Piola never accepted a university chair: however his pupil Francesco
Brioschi was the founder of the Politecnico di Milano. Brioschi mentored
Enrico Betti and Eugenio Beltrami. Ulisse Dini was pupil of Enrico Betti,
being also his successor as the chair of Mathematical Analysis and Geometry
at the Universit\`{a} di Pisa. Gregorio Ricci Curbastro was pupil of Ulisse
Dini, Eugenio Beltrami and Enrico Betti. Tullio Levi-Civita was pupil of
Gregorio Ricci Curbastro.

The strength of the Italian School of Mathematical Physics, Mathematical
Analysis and Differential Geometry has been weakened by two processes, one
which it shares with all other National Schools and in general with all
groups of scientists, the second one which is more peculiar to the Italian
Nation.

\begin{enumerate}
\item It happens very often that some theories need to be rediscovered and
reformulated several times in different circles before becoming a
universally recognized part of knowledge. For instance, the basic ideas of
functional analysis and its founding concept of functional (which goes back
to the calculus of variations and that can be defined with the sentence "a
function whose argument is a function") were already treated in the papers
by Erik Ivar Fredholm and in Hadamard's 1910 textbook and had previously
been introduced in 1887 by the Italian mathematician and physicist Vito
Volterra. The theory of nonlinear functionals was continued by students of
Hadamard, in particular Fr\'{e}chet and L\'{e}vy. Hadamard also founded the
modern school of linear functional analysis, further developed by Riesz and
the group of Polish mathematicians around Stefan Banach. However, Heisenberg
and Dirac did need to rediscover many parts of a theory already known and
they developed such a theory until the moment at which von Neumann could
recognize that actually Quantum Mechanics had been formulated in terms of
what he called Hilbert Spaces. Simple laziness or the difficulty of
understanding the formalism introduced by other authors, lack of time or of
economical means. All of these may lead some very brilliant scientists to
ignore results obtained by other scientists, which are nevertheless relevant
for their work. Many of the mathematicians listed in the previous genealogy
rediscovered many times the results which their predecessors had already
obtained because of the aforementioned first process. Such a process could
be called removal and/or ignorance of the results which appear not to be
relevant. This first removal process is indeed observed very often in the
history of science applied to the most various groups, independently of
their nationality\footnote{%
The authors are indebted to Prof. Mario Pulvirenti for having attracted
their attention to this first process and also for having recalled to them
the example concerning functional analysis.}, and the case of the
rediscovery of functional analysis is a striking example of its occurrence.

\item Napoleon favoured the birth of an Italian mathematical school, and
among many other actions he pushed Vincenzo Brunacci to write the first
Italian textbook in Mathematical Analysis. However he could not enforce in
the Italian School the habit, always followed by the French School, which
leads all French Scientists to recognize, to develop and to glorify the
contributions of their compatriots. Instead the Italian scientists always
preferred to follow the tradition of their predecessors, i.e. the scientists
of Greek language who developped the Hellenistic science (see \cite%
{russo2003}). Hellenic tradition is based on the intentional removal and
contempt of the contribution due to the compatriots and on the continuous
preference for the approval of foreign scientists. The described process
leads the members of a national group to consider the other national groups
always stronger, more qualified and more productive, while actively acting
to impeach the cultural, political and academic growth of the compatriot
scientists.
\end{enumerate}

The momentum given to the Italian School by Napoleon eventually lead to the
birth of Tensor Calculus, but was exhausted by the typical Italian negative
attitude towards compatriots, which was exemplified by the removal of
Levi-Civita, due to Mussolini's racial persecutions, from his chair in Rome,
that was immediately occupied by Signorini. Finally, however, it has to be
recognized at least that

i) the strict relationship between differential geometry and continuum
mechanics has been discovered and developed by the Italian School started by
Piola and culminating in Levi-Civita

ii) the great advancement of Riemannian geometry produced by the recognition
of the unicity of the parallel transport compatible with a Riemannian metric
(the so-called Levi-Civita Theorem) has its deep roots already in Piola's
works (recall the well-known concept of Piola's Transformation)

iii) Ulisse Dini's Theorem clarifies mathematically the concept of
constraint intensively used in the works of Piola. Indeed the crucial
concept of independent constraints (defined as those having non-singular
Jacobian) was clearified by Dini several decades after Piola had proven its
importance in continuum mechanics.

\end{document}